\documentclass [12 pt, leqno] {amsart}

\usepackage {amssymb}
\usepackage {amsthm}
\usepackage {amscd}

\usepackage[all]{xy}
\usepackage {hyperref}

\begin {document}

\theoremstyle{plain}
\newtheorem{proposition}[subsection]{Proposition}
\newtheorem{lemma}[subsection]{Lemma}
\newtheorem{korollar}[subsection]{Corollary}
\newtheorem{thm}[subsection]{Theorem}
\newtheorem*{thm*}{Theorem}

\theoremstyle{definition}
\newtheorem{definition}[subsection]{Definition}
\newtheorem{notation}[subsection]{Notation}

\theoremstyle{remark}
\newtheorem{bemerkung}[subsection]{Remark}

\numberwithin{equation}{subsection}

\newcommand{\eq}[2]{\begin{equation}\label{#1}  #2 \end{equation}}

\newcommand{\ml}[2]{\begin{multline}\label{#1}  #2 \end{multline}}

\newcommand{\nnal}[1]{\begin{align*} #1 \end{align*}}
\newcommand{\nneq}[1]{\begin{equation} \nonumber #1 \end{equation}}

\newcommand{\nnml}[1]{\begin{multline}\nonumber #1 \end{multline}}

%xypic

\newcommand{\arir}{\ar@{^{(}->}}
\newcommand{\aril}{\ar@{_{(}->}}

\newcommand{\are}{\ar@{>>}}

% Pfeile xr , xl 
\newcommand{\xr}[1] {\xrightarrow{#1}}
\newcommand{\xl}[1] {\xleftarrow{#1}}

% mathfrac, mathcal 

\newcommand{\mc}[1]{\mathcal{#1}}

% rm - Abkuerzungen 

\newcommand{\codim}{{\rm codim}\,}
\newcommand{\cd}{{\rm cd}}
\newcommand{\cone} {{\rm cone}}

\newcommand{\Hom} {{\rm Hom}}

% Abkuerzungen zu Bourbaki-Notation

% Zahlen
\newcommand{\Z} {\mathbb{Z}}

\newcommand{\C} {\mathbb{C}}

% affine und projektive Raueme
\newcommand{\A}{\mathbb{A}}

\newcommand{\G}{\mathbb{G}}

% mathcal Abkuerzungen 

% Sonst
\newcommand{\abs}[1]{\lvert#1\rvert}

\newcommand{\Zt}{\Z_{tr}}
\newcommand{\DMet}{{\bf DM}^{-}_{\text{\'et}}}
\renewcommand{\L}{\Lambda}
\newcommand{\Ltr}{\Lambda_{tr}}
\newcommand{\Het}[1]{H^{#1}_{\text{\'et}}}
\newcommand{\Aut}{{\rm Aut}}

\title[Cohomology of a linear subspace arrangement]{\'Etale cohomology of the complement of a linear subspace arrangement}
\author{Andre Chatzistamatiou}
\address{ University Duisburg-Essen \\ Fachbereich Mathematik Campus Essen \\ 45117 Essen \\ Germany}
\email{a.chatzistamatiou@uni-due.de}

\begin{abstract}
We prove a formula for the cup product on the $\ell$-adic cohomology of the complement of a linear subspace arrangement. 
\end{abstract}

\maketitle

\section*{Introduction}

Let $\mc{A}$ be a finite set of proper linear subspaces of the  $n$-dimensional affine space $V=\A_k^n$. 
We assume that $\mc{A}$ is closed under intersection and call it a linear subspace arrangement. In this paper we study
the $\ell$-adic \'etale cohomology ring of the complement $X=V-\cup_{A\in \mc{A}}A$ over an algebraically closed field $k$.  

By a sheaf theoretic approach Deligne, Goresky and MacPherson prove a decomposition 
\eq{zerlegungeinleitung}
{
H^*(X) \cong H^0(X)\oplus \bigoplus_{A\in \mc{A}} h^*(A)
} 
of the singular cohomology in the case $k=\C$, and of the $\ell$-adic \'etale cohomology if $k$ is
algebraically closed \cite[Exemple~1.11]{DGM}. The groups $h^*(A)$ can be described by combinatorial data as follows. 

We may consider $\mc{A}^+:=\mc{A}\cup \{ V \}$ as a category, with a unique arrow $A_1\xr{} A_2$
if $A_2$  contains $A_1$. If $A\in \mc{A}$ we denote by $\rangle A, V \langle$ the classifying space of the full subcategory 
$\{A'\in \mc{A}; A \subsetneq A'\}$ of $\mc{A}^+$. It is proven in \cite[1.10.1]{DGM} that  
\eq{kombinatorikeinleitung}
{
h^p(A) \cong H_{2 \cdot \codim A - p - 2} \left(\rangle A, V \langle;\Lambda \right) \quad \text{for all $p$,}
}  
where the right hand side means homology in the coefficients $\Lambda$ of \ref{zerlegungeinleitung}.  

In the case of singular cohomology the cup product is computed in \cite[Section~4]{DGM}. In particular, it is known that 
$$
h^*(A)\otimes h^*(B) \xr{\cup} h^*(C)
$$
vanishes if $A\cap B\not=C$ or if the intersection of $A$ and $B$ is not transverse. 
The proof uses singular chains in order to study the dual of the restriction to the diagonal  
\small$X \xr{} X \times X,$ \normalsize
and cannot be applied to the \'etale cohomology \cite[Remarque~4.3]{DGM}.
Restriction to the diagonal can be considered as a special case of a restriction map $H^*(X)\xr{} H^*(X')$ where 
$X'=X\cap V'$ is the intersection with a linear subspace $V'$ of $V$.
The decomposition \ref{zerlegungeinleitung} for singular cohomology has the following property 
\cite[Lemma~3.3]{DGM}: 
\begin{itemize}
\item [(P)] 
For all $M\in \mc{A}$ the composition $h^*(M) \xr{} H^*(X) \xr{} H^*(X')$ vanishes if $M\cap V'=\emptyset$  or if 
the intersection is not transverse.
\end{itemize}

In this paper we construct another decomposition of the \'etale cohomology $H^*(X)$ which has the same 
combinatorial description \ref{kombinatorikeinleitung} and satisfies the condition (P).
For this decomposition we can prove the conjectured formulas of \cite{DGM} for the cup product. If the decomposition from 
\cite{DGM} satisfies (P) then both decompositions agree. Unfortunately we are not able to prove this.   

In order to construct our decomposition we use Voevodsky's work on the category of complexes of \'etale sheaves with transfers, 
and the localization $\DMet$ where $\A^1$-weak equivalences become an isomorphism. For finite coefficients $\Lambda = \Z/n$ and 
$n^{-1}\in k$ we have \cite[Prop.~10.7]{Voevodskylecture}
$$
\Hom_{\DMet}(\Zt(X),\L[p]) \cong \Het{p}(X,\L) \quad \text{for all $p$}.
$$ 
We find a resolution of the sheaf $\Zt(X)$ and decompose this complex modulo cones of $\A^1$-weak equivalences (Section~\ref{sec1},\ref{sec2}). 
In Section \ref{sec3} we study the restriction to linear subspaces, and in Section \ref{sec4} we prove the formula for the cup product.

We have to introduce some notation to state the formula for the cup product. We denote by $\langle A,B \rangle$ 
(resp. $\rangle A,B \rangle$, $\langle A,B \langle$) the classifying space of \small$\{ A'; A \subset A' \subset B \}$ \normalsize 
(resp. \small$\{ A'; A \subsetneq A' \subset B \}$, $\{ A'; A \subset A' \subsetneq B \}$\normalsize). 
Let $\mc{A}^+\times \mc{A}^+$ be the product category; the functor $(A'_1,A'_2)\xr{}A'_1\cap A'_2$ 
induces a morphism \small$\tau\colon\langle A_1\times A_2,V\times V \rangle \xr{} \langle A_1\cap A_2,V \rangle$ \normalsize for $A_1,A_2\in \mc{A}^+$.

\begin{thm*}
Let $A_1,A_2\in \mc{A}^+, \L=\Z/n$ and $n^{-1}\in k.$ We assume that $k$ is algebraically closed, 
%%%%%>>>>>>>>>> revision 1
and that the arrangement $\mc{A}$ is defined over a subfield $k_0\subset k$.
\begin{itemize}
\item[(i)] We have a decomposition of $Gal(k/k_0)$-modules
\small$$
\Het{p}(X,\L)= \bigoplus_{A\in \mc{A}^+-\{\emptyset\}} 
H_{2\cd A-p}(\langle A,V \rangle,\rangle A,V \rangle\cup \langle A,V \langle;\L)\otimes \mu_n^{\otimes -\cd A} 
$$ \normalsize
for all $p$. Here $\cd$ denotes the codimension in $V$, and the action of $Gal(k/k_0)$ on
$H_{2\cd A-p}(\langle A,V \rangle,\rangle A,V \rangle\cup \langle A,V \langle;\L)$ is trivial.  
%%%%%<<<<<<<<<< revision 1
\item[(ii)] The cup product on two summands from (i) for $A_1,A_2$ vanishes if $A_1\cap A_2=\emptyset$ or 
$\cd (A_1\cap A_2) < \cd A_1 + \cd A_2$. Otherwise it is (up to twist with $\mu_n^{\otimes -\cd A_1-\cd A_2}$) the composition: 
\begin{footnotesize}
$$
\begin{CD}
H_{2\cd A_1-p_1}(\langle A_1,V \rangle,\rangle A_1,V\rangle \cup \langle A_1,V \langle ;\L) \otimes  
H_{2\cd A_2-p_2}(\langle A_2,V \rangle , \rangle A_2,V \rangle \cup \langle A_2,V \langle ;\L) \\
@VV\text{external product}V\\
H_{2\cd A_1 \cap A_2-p_1-p_2}(\langle A_1\times A_2,V^2 \rangle , \rangle A_1\times A_2,V^2 \rangle \cup \langle A_1\times A_2,V^2 \langle ;\L)  \\
@VV{H_*(\tau)}V\\
H_{2\cd A_1\cap A_2-p_1-p_2}(\langle A_1\cap A_2,V \rangle , \rangle A_1\cap A_2,V \rangle \cup \langle A_1\cap A_2,V \langle ;\L).
\end{CD}
$$
\end{footnotesize}
%% \begin{scriptsize}
%% $$
%% \xymatrix @-1pc
%% {
%% H_{2\cd A_1-p_1}([A_1,V],]A_1,V]\cup [A_1,V[;\L) \otimes  H_{2\cd A_2-p_2}([A_2,V],]A_2,V]\cup [A_2,V[;\L) \ar[d]
%% \\
%% H_{2\cd A_1 \cap A_2-p_1-p_2}([A_1\times A_2,V^2],]A_1\times A_2,V^2]\cup [A_1\times A_2,V^2[;\L) \ar[d]^{\tau}
%% \\
%% H_{2\cd A_1\cap A_2-p_1-p_2}([A_1\cap A_2,V],]A_1\cap A_2,V]\cup [A_1\times A_2,V[;\L).
%% }$$ \end{scriptsize}
\end{itemize} 
\end{thm*}   

\subsection*{Acknowledgments}
I thank H. Esnault for introducing me to this subject and for her interest in my work. I am grateful to Y. Andr\'e for valuable discussions. 
This paper is written during a stay at the Ecole normale sup\'erieure which is supported by a fellowship within the Post-Doc program of the Deutsche Forschungsgemeinschaft (DFG).
I thank the Ecole normale sup\'erieure for its hospitality.
  
\section{Resolution of the complement}
\label{sec1} 
\subsection{}
Let $k$ be a field. We work with $Sh_{Nis}(Cor_k)$, the category of Nisnevich sheaves with transfers (see 
\cite{Voevodskylecture}, Lecture 12). All results will be true in the \'etale topology as well.
We denote by $\Z_{tr}(V)$, for a smooth $k$-scheme $V$, the  Nisnevich sheaf
with transfers defined by $\Z_{tr}(V)(Y)= Cor(Y,V)$.

\subsection{} \label{notation}
Let $V$ be a smooth $k$-scheme and $\mc{A}$ a set of closed and proper 
subsets of $V$. We assume that intersections of elements in $\mc{A}$ are
contained in $\mc{A}$, i.e. $A_1,A_2\in \mc{A}\Rightarrow A_1\cap A_2\in \mc{A}$, 
and set $\mc{A}^+:=\mc{A} \cup\{ V \}$.
We may view $\mc{A}^+$ as a category, with an arrow $A_1\rightarrow A_2$ if $A_1\subset A_2$. 
We denote by
$$
m \mapsto \mc{A}^+(m):=\{ A_0\rightarrow A_1\rightarrow \dots \rightarrow A_m ; A_i \in \mc{A}^+\}
$$
the associated simplicial set, i.e. the nerve of $\mc{A}^+$, with face operators 
\small
$$
\partial_{i}(A_0\rightarrow \dots \rightarrow A_m)=
\begin{cases} 
A_1\rightarrow \dots \rightarrow A_m   &\text{if $i=0$,} \\ 
A_0\rightarrow \dots  A_{i-1} \rightarrow A_{i+1} \dots \rightarrow A_m  &\text{if $0<i<m$,} \\
A_0\rightarrow \dots \rightarrow A_{m-1}  &\text{if $i=m$.}  
\end{cases}
$$\normalsize 
Let $\mc{A}^+(m)_{\text{nd}} \subset \mc{A}^+(m)$ be the subset of non-degenerated elements, i.e. 
$A_i\not=A_{i+1}$ for all $i$.

\subsection{} 
The complex $S(V,\mc{A})$ is defined by 
$$
S(V,\mc{A})^{m}:= \bigoplus_{A_0\rightarrow \dots \rightarrow V \in \mc{A}^+(m)_{\text{nd}}} \Z_{tr}(V)/\Z_{tr}(V-A_0)
$$
with differential $d=\sum_{i} (-1)^{i} \partial^{i}$, where $ \partial^{i} $ is defined by 
$$
{\rm pr}_A \circ \partial^i :=  (  \Z_{tr}(V)/\Z_{tr}(V-\partial_i(A)_0) \xr{{\rm Q}}  \Z_{tr}(V)/\Z_{tr}(V-A_0)  )   \circ {\rm pr}_{\partial_i(A)}   
$$
for all $A \in \mc{A}^+(m)_{\text{nd}}$. Here, ${\rm pr}_A$ denotes the projection to the direct summand for $A$, and $Q$ denotes the quotient-morphism. 

We have a morphism of complexes
\eq{morphismUVA}
{
\Z_{tr}(V-\cup_{A\in \mc{A}} A) \xr{} S(V,\mc{A}) 
}
induced by $\Z_{tr}(V-\cup_{A\in \mc{A}} A)\subset \Z_{tr}(V)=S(V,\mc{A})^{0}$.

\begin{lemma}\label{lemma:aufloesung}
The morphism \ref{morphismUVA} is a quasi-isomorphism.
\begin{proof}
Clearly, we have $H^0(S(V,\mc{A}))=\Z_{tr}(V-\cup_{A\in \mc{A}} A)$; we show that $H^i(S(V,\mc{A}))=0$ for $i>0$.

Let $Y$ be the spectrum of a local henselian ring. The cohomology $H^i(S(V,\mc{A})(Y))$ is generated by elements of the form  
\eq{zutrivialisierendeKlasse}
{
\alpha=\sum_{A \in \mc{A}^+(i)_{\text{nd}}, A_{i}=V} c_{A} Z,
}
where $c_{A}\in \Z$ and $Z\in \Z_{tr}(V)(Y)$ is an irreducible cycle such that 
$$
\sum_{j=0}^{i} (-1)^{j} c_{\partial_jB} = 0
$$
for all $B\in \mc{A}^+(i+1)_{\text{nd}}$ with $Z \cap (Y\times B_0)\neq \emptyset$. 
Every intersection of elements in the set $\mc{A}_Z := \{A \in \mc{A}; Z \cap (Y\times A)\neq \emptyset  \}$ is contained in $\mc{A}_Z$ again.
Indeed, if $s\in Z$ denotes the closed point, then  
$$
Z \cap (Y\times A)\neq \emptyset \Leftrightarrow s \in Y\times A.
$$  
In particular, $\mc{A}_Z$ contains a minimal element $M$. Let $\abs{\mc{A}_Z}$ be the classifying space of the category $\mc{A}_Z$. The complex
$$
C^dS = \bigoplus_{A\in \mc{A}_Z(d)_{\text{nd}}} \Z
$$ 
computes the singular cohomology of $\abs{\mc{A}_Z}$, and using the cycle $Z$ we get a morphism 
%%%>>>>>>>>>  revision 1
$$
{\rm cone}(\Z \xr{} C^*S)[-1] \xr{\phi_Z} S(V,\mc{A})(Y). 
$$
Explicitly, the map $\phi_Z$ is 
$$
\Z \xr{} S(V,\mc{A})(Y)=\Z_{tr}(V)(Y); \quad 1 \mapsto Z,
$$ 
in degree $=0$, and in degree $d>0$, $\phi_Z^{d}:C^{d-1}S \xr{} S(V,\mc{A})^{d}(Y)$ is 
defined by 
$$
{\rm pr}_A \circ \phi_Z^{d} = \begin{cases} 0 & \text{if $\partial_d A\not\in \mc{A}_Z(d-1)_{{\rm nd}}$} 
\\ I^A_Z \circ {\rm pr_{\partial_d A}} & \text{otherwise,}  \end{cases}
$$ 
for $A\in \mc{A}^+(d)_{{\rm nd}}$ such that $A_d=V$. Here the map 
$$I^A_Z: \Z \xr{} \Z_{tr}(V)(Y)/\Z_{tr}(V-A_0)(Y)$$ 
maps $1$ to the cycle $(-1)^d Z$. By definition of $\mc{A}_Z$ the image of $\phi_Z$ contains the
elements of \ref{zutrivialisierendeKlasse}.
%%The shift [-1] was missing in previous version; now ${\rm cone}(\Z \xr{} C^*S)[-1]$ starts 
%%with $\Z$ in degree =0$. 
%%Note that $H^0(C^*S)\xr{} H^1(S(V,\mc{A})(Y))$, so cone construction is needed. Otherwise  
%%one has to consider $(Z)_{A_0 \xr{} V}\in H^1(S(V,\mc{A})(Y))$ (this is the image of a
%%generator of $H^0(C^*S)$) separately, which is not difficult: it is the image of 
%%$Z\in S(V,\mc{A})^0(Y)$ by the differential. 
%%%<<<<<<<<<  revision 1

Since $\abs{\mc{A}_Z}$ is a cone over the minimal element $M$, the reduced singular cohomology vanishes and  
$
{\rm cone}(\Z \xr{} C^*S)
$
is acyclic. Therefore, the class of the element \ref{zutrivialisierendeKlasse} is trivial.
\end{proof}
\end{lemma}

\subsection{} \label{funkA}
Let $V'$ be another smooth $k$-scheme and $\mc{A}'$ a finite set of closed, proper subsets of $V'$. 
Furthermore, let
\begin{enumerate}
\item $\phi\colon \mc{A} \xr{} \mc{A}'$ be a map preserving inclusions,
\item $f\in \Z_{tr}(V)(V')=Cor(V',V)$ be a finite correspondence.
\end{enumerate}
We assume that $f$ is given by $f=\sum_{Z} c_Z Z$, where $c_Z\in \Z$ and \linebreak $Z\subset V\times V'$ is irreducible, and 
such that:
$$
Z\cap (V'  \times A ) \subset \phi(A) \times V \quad \text{for all $c_Z\neq 0$ and $A\in \mc{A}$.}
$$
Then we get a morphism of complexes 
\eq{funkR}
{
S(f,\phi)\colon  S(V',\mc{A}') \xr{} S(V,\mc{A}) 
}
defined by  
$$
{\rm pr}_{A} \circ S(f,\phi) = 
\begin{cases} 
\bar{f} \circ {\rm pr}_{\phi(A)} &  \text{if $\phi(A)\in {\mc{A}'}^+(m)_{\text{nd}}$} \\
0 & \text{if $\phi(A)$ is degenerated}
\end{cases} 
$$
for all $A\in \mc{A}^+(m)_{\text{nd}}$. 
Here  
$$
\bar{f}\colon \Z_{tr}(V')/\Z_{tr}(V'-\phi(A)_0) \xr{} \Z_{tr}(V)/\Z_{tr}(V-A_0) 
$$
is the morphism induced by $f$, and we extend the map $\phi$ by 
$$
\phi(A_0\rightarrow \dots \rightarrow A_{m-1}\rightarrow V):= \phi(A_0)\rightarrow \dots \rightarrow \phi(A_{m-1})\rightarrow V'.
$$

Obviously the following diagram is commutative: 
\eq{nat}
{
\xymatrix
{
S(V',\mc{A}') \ar[r]^{S(f,\phi)} 
&
S(V,\mc{A})  
\\
\Zt(V'-\cup_{A'\in \mc{A}'} A') \ar[u]^{{\rm quis}} \ar[r]^{f}
&
\Zt(V-\cup_{A\in \mc{A}} A). \ar[u]_{{\rm quis}}
}
}

\section{Motivic decomposition of the complement of a subspace arrangement} 
\label{sec2}

\subsection{} \label{betanotation}
Let $V=\A^n_k$ be the affine space and $\mc{A}$ a finite set of proper, linear subspaces in $V$ (we allow $\emptyset$
to be a linear subspace). As in Section \ref{notation} we assume that intersections of elements of $\mc{A}$ are contained 
in $\mc{A}$.

The construction of our motivic decomposition depends on a suitable transverse arrangement for $\mc{A}$. 
Let $\mc{T}$ be a set of linear subspaces in $V$ and $\beta\colon \mc{A}^{+} \xr{} \mc{T}$ a bijection. 
We assume that the following conditions for $\beta$ are satisfied:
\begin{enumerate}
\item[(T1)] \label{beta1} $\beta(\emptyset)=V$ if $\emptyset \in \mc{A}$,   
\item[(T2)] \label{beta2} $A_1\subset A_2 \Leftrightarrow  \beta(A_1) \supset \beta(A_2),$
\item[(T3)] \label{beta3} For all $ A\in \mc{A}^{+} - \{\emptyset \}$ we have: ${\rm codim}(A)={\rm dim}(\beta(A))$ and $\beta(A)\cap A$ is a point. 
\item[(T4)] \label{beta4} For all $A,B \in \mc{A}^+$ we have: $\beta(A)\cap A\cap B\not= \emptyset \Rightarrow A\subset B$.
\end{enumerate} 
By (T3) the subspaces $A$ and $\beta(A)$ intersect transversely in one point. Condition (T4) implies that $\beta(A)$ intersects transversely with  all $B\in \mc{A}$, that contain that point.
One way to construct a $\beta$ is as follows. 
Choose a point $x\in V(k)-\cup_{A\in \mc{A}} A(k)$ and consider $V$ as a vector space with $x$ as zero. For a non-degenerated bilinear form $B$ we set 
$$
 \beta_B(A):=\{w\in V ; B(w,v)=0 \; \text{for all $v\in T(A)$} \},
$$ 
where $T(A)$ is the tangent space of $A$, i.e. the linear subspace parallel to $A$ that contains $x$. 
It is easy to see that condition (T3) and (T4) for $\beta_B$ are satisfied for an open, non-empty set
of bilinear forms. Therefore, transverse arrangements exist if $k$ is a field with infinitely many elements 
(since non-degenerated bilinear forms are an open set in an affine space). If $k$ is finite then  
one finds transverse arrangements after a base change to a finite extension of $k$. 

\subsection{}\label{definitionSbeta}
If $\beta$ is a transverse arrangement for $\mc{A}$, we denote by $S_{\beta}(V,\mc{A})$ the complex  
$$
S_{\beta}(V,\mc{A})^{m}:= \bigoplus_{A \in \mc{A}^+(m)_{\text{nd}}, A_{m}=V} \Z_{tr}(\beta(A_0))/\Z_{tr}(\beta(A_0)-A_0)
$$
with differential $d=\sum_{i>0} (-1)^{i} \partial^{i}$, where $\partial^{i}$ is defined as in Section \ref{notation} 
(here we have $\partial^0=0$).
The inclusions $\beta(A)\subset V$ induce a morphism
\eq{SbetaS}
{
S_{\beta}(V,\mc{A}) \xr{} S(V,\mc{A}). 
}
Here we used $\beta(B)\cap A=\emptyset$ for all $A\subsetneq B$ (which follows from (T4)). 

\subsection{}\label{seckompbb'}
Let $V'$ be another  affine space and $\mc{A}'$ a set of linear subspaces in $V'$ as in \ref{betanotation}.
Furthermore, let $f\colon V'\xr{} V$ be an affine-linear map and $\phi\colon\mc{A} \xr{} \mc{A}'$ as in \ref{funkA}, i.e. 
for all $A\in \mc{A}$ we have $f^{-1}(A)\subset \phi(A)$. If $\beta'$ is a transverse arrangement for $\mc{A}'$ 
that satisfies  
\eq{kompbb'}
{
f(\beta'(\phi(A))) \subset \beta(A)
} 
for all $A\in \mc{A}^+$, then we get a commutative diagram  
\eq{diakompbb'}
{
\xymatrix
{
S_{\beta'}(V',\mc{A}') \ar[r] \ar[d] 
&
S_{\beta}(V,\mc{A}) \ar[d]
\\
S(V',\mc{A}') \ar[r]^{S(f,\phi)} 
&
S(V,\mc{A}). 
}
}

\begin{bemerkung}
In general, for given arrangements $\mc{A}',\mc{A},$ and given affine linear 
map $f\colon V' \xr{} V$ there are no transverse arrangements $\beta',\beta$, such that
\ref{kompbb'} holds. Indeed, let $V=\A^2$ and $\mc{A}=\{L_1,L_2,s\}$ where $L_1,L_2$ are two lines that 
intersect in the point $s$. Let $V'$ be another line through $s$, and $\mc{A}'=\{s\}$. For every
$\beta'$ we have $\beta'(s)=V'$ and $V'\not=\beta(L_1)$ since $s\in V'$.
   
\end{bemerkung}

\subsection{}
Let $\mathbf{D}^-=\mathbf{D}^{-}(Sh_{Nis}(Cor_k))$ be the derived category of bounded above 
complexes in $Sh_{Nis}(Cor_k)$.
We recall Voevodsky's definition of $\A^1$-weak equivalence \cite[Definition~9.2]{Voevodskylecture}. 
A morphism $f$ in $\mathbf{D}^-$ is called an $\A^1$-weak equivalence if the cone of $f$ is in 
the smallest thick subcategory $\mc{E}_{\A}$ such that:
\begin{enumerate}
\item $\cone(\Z_{tr}(V\times \A^1) \xr{} \Z_{tr}(V))$ is in $\mc{E}_{\A}$ for all smooth scheme $V$,  
\item $\mc{E}_{\A}$ is closed under any direct sum that exists in $\mathbf{D}^-$. 
\end{enumerate} 
The $\A^1$-weak equivalences form a saturated, multiplicative system, which we denote by $W_{\A}$. 
The localization $\mathbf{D}^-[W_{\A}^{-1}]$ is  a triangulated category 
and a morphism $f$ in $\mathbf{D}^-$ becomes an isomorphism in $\mathbf{D}^-[W_{\A}^{-1}]$ if and only if $f\in W_{\A}$.

\begin{lemma} \label{lemma:iso}
The morphism  
$$
S_{\beta}(V,\mc{A}) \xr{} S(V,\mc{A}) 
$$
from \ref{SbetaS} is an $\A^1$-weak equivalence.
\begin{proof}
It is enough to prove that  
\eq{claimlemmaiso}
{\Z_{tr}(\beta(A))/\Z_{tr}(\beta(A)-A) \xr{} \Z_{tr}(V)/\Z_{tr}(V-A)}
is an $\A^1$-weak equivalence for all $A\in \mc{A}^+ - \{\emptyset\}$.

For $A\in \mc{A}^+ - \{\emptyset\}$ let $\epsilon_A\colon V\rightarrow \beta(A) \times A$ be the affine linear map that 
identifies $A$ with $x_A \times A$ and 
$\beta(A)$ with $\beta(A) \times x_A$, where $x_A=\beta(A)\cap A$ is the intersection point.   
The morphism $\Z_{tr}(\beta(A)) \xr{} \Z_{tr}(V)$ (resp. $\Z_{tr}(\beta(A)-A) \xr{} \Z_{tr}(V-A)$) 
is an $\A^1$-weak equivalence, because it is the inverse (in $\mathbf{D}^-[W_{\A}^{-1}]$) of the  
projection $ pr_1\circ \epsilon_A\colon \Z_{tr}(V)\xr{} \Z_{tr}(\beta(A)) $  
(resp. $ pr_1\circ \epsilon_A\colon \Z_{tr}(V-A)\xr{} \Z_{tr}(\beta(A)-A) $). 
The claim follows from the 5-lemma.
\end{proof}
\end{lemma}

\subsection{} \label{arrangementopology}
Working with arrangements it is convenient to introduce a topology on $\mc{A}^+$. 

%%%>>>>>>>>> revision 1
More generally, let $\mc{C}$ be a category such that the objects form a set 
(by abuse of notation we write ${\rm Ob}(\mc{C})=\mc{C}$).  
We define a topology on the set of objects by setting: a set $U\subset \mc{C}$ is open
if for every $B\in U$ and $C\in \mc{C}$ the existence of an arrow $B\xr{} C$ implies that
$C\in U$.

Closed sets $Y$ satisfy: every morphism with target in $Y$ has its source in $Y$. 

The smallest open set containing $M\in \mc{C}$ is 
$\{N\in \mc{C}; \exists M\xr{}N\}$ (in the case $\mc{C}=\mc{A}^+$ we denote 
this set by $[M,V]$). Similarly, the smallest closed set
containing $N\in \mc{C}$ is $\{M\in \mc{C}; \exists M\xr{}N\}$.

If $\mc{C}'\subset \mc{C}$ is a locally closed set, i.e. $\mc{C}'=U\cap Y$ for an open set
$U$ and a closed set $Y$, and $M\xr{} A \xr{} N$ is the composite of two arrows in 
$\mc{C}$ with $M, N \in \mc{C}'$, then $A\in \mc{C}'$. 
Indeed, $N\in Y$ implies $A\in Y$, and $M\in U$ implies $A\in U$.   

For a full subcategory $\mc{C}'\subset \mc{C}$ the topology induced by $\mc{C}$ equals the
topology of $\mc{C}'$.
%%%<<<<<<<<< revision 1

%%%%>>>>>>>>> revision 1
\begin{notation} \label{notationclassifying}
We consider subsets $\mc{C}\subset \mc{A}^+$ as full subcategories of $\mc{A}^+$. We denote
by 
$$
m \mapsto \mc{C}(m):=\{ C_0\rightarrow C_1\rightarrow \dots \rightarrow C_m ; C_i \in \mc{C}\}
$$
the nerve of $\mc{C}$ and by $\mc{C}(m)_{\rm nd}$ the non-degenerated elements.
We write $\abs{\mc{C}}$ for the classifying space; it is a closed subset of $\abs{\mc{A}^+}$. 
\end{notation}

For $\mc{C}\subset \mc{A}^+$ the complex 
\eq{SC}
{
S^{\mc{C}}_{m}:= \bigoplus_{C \in \mc{C}(m)_{\text{nd}}} \Z [C],
}
with differential $\sum_{i} (-1)^{i} \partial_{i}$, computes the singular homology 
of $\abs{\mc{C}}$.
%%%<<<<<<<<<< revision 1

\subsection{}
Given two locally closed subsets $\mc{M},\mc{N}\subset \mc{A}^+$ 
we define $S^{\mc{M},\mc{N}}$ to be the complex 
\eq{SMN}
{
S^{\mc{M},\mc{N}}_{m}:= \bigoplus_{\substack{A \in \mc{A}^+(m)_{\text{nd}} \\ A_0\in \mc{M}, A_m\in \mc{N}}} \Z [A]
}
with differential $d=\sum_{i} (-1)^{i} \partial_{i}$, where $\partial_0(A)=0$ if $A_{1}\not\in \mc{M}$, 
and $\partial_m(A)=0$ if $A_{m-1}\not\in \mc{N}$. 

%%%%>>>>>>>>>>>>> revision 1
Writing $\mc{M}=U\cap Z$ and $\mc{N}=W\cap Y$,
with open sets $U,W$ and closed sets $Z,Y$, and setting $\mc{C}=U\cap Y$, 
we see that $S^{\mc{M},\mc{N}}$ is a quotient 
$$
S^{\mc{C}} \xr{} S^{\mc{M},\mc{N}}.
$$
The kernel is generated by chains $C\in \mc{C}(m)_{{\rm nd}}$ such that 
$C_0\not\in Z$ or $C_m\not\in W$. Since $Z$ is closed in $\mc{A}^+$, the set
$\mc{C}-Z=\mc{C}-\mc{M}$ is open in $\mc{C}$; similarly $\mc{C}-W=\mc{C}-\mc{N}$
is closed in $\mc{C}$. Therefore, $C_0\not\in Z$ implies $C\in (\mc{C}-\mc{M})(m)_{{\rm nd}}$,
and $C_m\not\in W$ implies $C\in (\mc{C}-\mc{N})(m)_{{\rm nd}}$. We get an exact sequence
\eq{SMNaufloesung}
{
0 \xr{} S^{\mc{C}-(\mc{M}\cup \mc{N})} \xr{} S^{\mc{C}-\mc{M}} \oplus S^{\mc{C}-\mc{N}} \xr{} S^{\mc{C}}  \xr{} S^{\mc{M},\mc{N}} \xr{} 0.
}  
From the sequence we see that 
$$
H_m(S^{\mc{M},\mc{N}}\otimes \Lambda) = H_m(\abs{\mc{C}},\abs{\mc{C}-\mc{M}} \cup \abs{\mc{C} -\mc{N}}; \Lambda),
$$
i.e.~$S^{\mc{M},\mc{N}}\otimes \Lambda$ computes the relative singular homology of $\abs{\mc{C}}$ 
relative to the closed set $\abs{\mc{C}-\mc{M}} \cup \abs{\mc{C} -\mc{N}}$ 
(with coefficients $\Lambda$).

In the case $\mc{M}=\{M\},\mc{N}=\{V\}$, we have $U=[M,V]$ and $Y=\mc{A}^+$, so that 
$\mc{C}=[M,V]$. 
The spaces $\abs{\mc{C}},\abs{\mc{C}-\mc{M}}$ and $\abs{\mc{C}-\mc{N}}$ are contractible 
if $V\neq M$ ($\mc{C}-\mc{M}$ has a maximal element $V$, and $\mc{C}-\mc{N}$ has a minimal 
element $M$). Thus
$$
H_{p}(\abs{\mc{C}},\abs{\mc{C}-\mc{M}}\cup \abs{\mc{C}-\mc{N}}) \cong \tilde{H}_{p-2}(\abs{\mc{C}-(\mc{M}\cup \mc{N})}) \quad \text{for all $p$.}
$$  
%%%%<<<<<<<<<<< revision 1

\subsection{} 
Let $\mc{M},\mc{N},\mc{P}$ be three locally closed subsets in $\mc{A}^+$; we have a morphisms  
$$
\Hom(S^{\mc{M},\mc{N}},\Z) \otimes S^{\mc{M},\mc{P}} \xr{} S^{\mc{N},\mc{P}} 
$$ 
($\Hom$ denotes the $\Hom$-complex) defined by  
\nnml
{
[A'_0\xr{} \dots \xr{} A'_{m_1}]^*  \otimes [A_0\rightarrow \dots \xr{} A_{m_2}] \mapsto \\
\begin{cases}
0 & \text{if $A'_i\neq A_i$ for some $i\leq m_1$,} \\
 [A_{m_1}\rightarrow \dots \xr{} A_{m_2}] & \text{otherwise.} 
\end{cases}
}
This morphism comes from the cap product for singular chains. 

Taking duals we get 
\eq{kom}
{
  \Hom(S^{\mc{N},\mc{P}},\Z) \xr{} S^{\mc{M},\mc{N}} \otimes \Hom(S^{\mc{M},\mc{P}},\Z). 
}

\subsection{} \label{Gm}
If $F$ is a sheaf, the Suslin complex $C_*(F)$ \cite[Definition~2.14]{Voevodskylecture} is defined
by $C_m(F)(U):=F(U\times \Delta^m)$, where $\Delta^m:=\{(t_0,\dots,t_m); \\ \sum_{i\geq 0} t_i = 0 \}$ is the
algebraic simplex. The differential is given by $\sum_i (-1)^{i} \partial_i$, where, by abuse of notation, 
$\partial_i$ denotes the restriction $C_m(F)\rightarrow C_{m-1}(F)$ to the face $t_i=0$. We have a natural morphism
\eq{natMorphFCF}
{
F \xr{} C_*(F).
}
Following Voevodsky, we have
$$
\Z(1):={\rm cone}(C_*(\Z)\xr{} C_*(\Z_{tr}(\G_m)))[-1]
$$
by using the point $1\in \G_m$ \cite[Definition~3.1]{Voevodskylecture}.
Since $F\xr{} C_*(F)$ and $\Z_{tr}(\A^1)\xr{} \Z$ are $\A^1$-weak equivalences \cite[Lemma~14.4]{Voevodskylecture}, we can write 
\nnal
{
\Z(1) &= {\rm cone}(\Z\xr{} \Z_{tr}(\G_m))[-1] \\
      &= {\rm cone}(\Z_{tr}(\G_m)\xr{} \Z_{tr}(\A^1))[-2]
}
in $\mathbf{D}^-[W_{\A}^{-1}]$.
We have $\Z(q):=\Z(1)^{\otimes^{tr}q}=$
\nnml
{
[\Z_{tr}(\G_m^{\times q}) \xr{} \oplus_i \Z_{tr}(\G_m^{\times i}\times \A^1 \times \G_m^{\times q-i-1}) \xr{} \dots \xr{} \\
   \oplus_i \Z_{tr}(\A^{i}\times \G_m \times \A^{q-i-1}) \xr{} \Z_{tr}(\A^{q})],
}
where $\Z_{tr}(\G_m^{\times q})$ is in degree $q$ and all arrows come from the inclusion $\G_m\subset \A^1$ (with adequate sign). 
The sets $\{(x_j) \in \A^{n}-\{ 0 \}; x_i \neq 0 \}$ cover $\A^n-\{ 0 \}$, and therefore 
\nnml
{
[\Z_{tr}(\G_m^{\times q}) \xr{} \oplus_i \Z_{tr}(\G_m^{\times i}\times \A^1 \times \G_m^{\times q-i-1}) \xr{} \dots \xr{} \\
   \oplus_i \Z_{tr}(\A^{i}\times \G_m \times \A^{q-i-1}) \xr{} \Z_{tr}(\A^{q}-\{0\})]
}
is exact in the Nisnevich (or \'etale) topology \cite[Proposition~6.12]{Voevodskylecture}. This implies 
$$
\Z(q)[2q]\cong \Z_{tr}(\A^{q})/\Z_{tr}(\A^{q}-\{0\})
$$ 
in $\mathbf{D}^-[W_{\A}^{-1}]$. Choosing a linear isomorphism $\epsilon\colon \beta(M)\xr{} \A^{\codim M}$ such that
$\epsilon(\beta(M)\cap M)=0$, we get  
\eq{Gysin?}
{
\Z(\codim M)[2\codim M]\cong \Z_{tr}(\beta(M))/\Z_{tr}(\beta(M)-M).
}

\begin{bemerkung} 
If $k$ is a perfect field there is a natural isomorphism as in \ref{Gysin?} from the
Gysin-sequence \cite[Proposition~3.5.4]{Voevodsky}. 
\end{bemerkung}

\subsection{} 
From the definition of $S_{\beta}(V,\mc{A})$ in \ref{definitionSbeta}, we can see immediately the decomposition:
\eq{ZerlegungSbeta}
{
S_{\beta}(V,\mc{A})=\bigoplus_{M\in \mc{A}^+}  \Z_{tr}(\beta(M))/\Z_{tr}(\beta(M)-M) \otimes  \Hom(S^{\{M\},\{V\}},\Z).
}

From \ref{ZerlegungSbeta}, Lemma \ref{lemma:aufloesung} and \ref{lemma:iso}, we get the following theorem. 
\begin{thm}\label{thm:zerlegung}
Let $V=\A^n_k$ be the affine space and $\mc{A}$ a finite set of proper linear subspaces in $V$, such that intersections of elements of $\mc{A}$ are contained in $\mc{A}$.
If $\beta \colon \mc{A}^{+} \xr{} \mc{T}$ is a transverse arrangement as in Section \ref{betanotation} then there is an 
isomorphism $I_{\beta}$ in $\mathbf{D}^-[W_{\A}^{-1}]$:
\small
\nnal
{
\Z_{tr}(V-\cup_{A\in \mc{A}}A) &\xr{I_{\beta}}  
 \bigoplus_{M\in \mc{A}^+}  \Z_{tr}(\beta(M))/\Z_{tr}(\beta(M)-M) \otimes  \Hom(S^{\{M\},\{V\}},\Z) \\
 &\cong \bigoplus_{M\in \mc{A}^+-\{\emptyset\}} \Z(\codim M)[2\codim M] \otimes  \Hom(S^{\{M\},\{V\}},\Z).
}\normalsize
\end{thm}

\begin{bemerkung}
In general, the isomorphism $I_{\beta}$ in the Theorem depends on the choice of the transverse arrangement.
\end{bemerkung}

%%%%>>>>>>>>>>>>> revision 1
\subsection{} \label{Galoisop}
If the arrangement is defined over a subfield $k_0\subset k$, then there is a natural 
``operation'' of the automorphism group $G=\Aut_{k_0}(k)$ on $\Z_{tr}(V-\cup_{A\in \mc{A}}A)$, which 
will be defined below. Essentially, this operation induces the operation of $G$ on \'etale 
cohomology. Since our decomposition of \'etale cohomology will be induced by $I_{\beta}$
we have to compute how $I_{\beta}$ transforms under $G$. 

For every smooth scheme $X$ over $k$ and $g\in G$ we form the fiber square
$$
\xymatrix
{
X_g \ar[r]^{g} \ar[d] &  X \ar[d] \\
k \ar[r]^{g} &                k.  
}
$$     
We note that $(X_g)_h=X_{g\circ h}$.

We define a homomorphism $G\xr{} \Aut(Cor_k)$ by 
\begin{equation*}
\begin{split}
X &\mapsto  X_{g^{-1}} \\
Cor(X,Y) & \xr{} Cor( X_{g^{-1}},  Y_{g^{-1}}), \; Z\mapsto  (g\times g)_* Z.  
\end{split}
\end{equation*}
Presheaves with transfers are functors from $Cor_k$ to abelian groups, so we 
get an homomorphism $G\xr{} \Aut(PreSh(Cor_k)), g\mapsto g_*,$ with $g_*F(Y):=F(Y_{g})$ 
for every presheaf $F$. It is easy to see that $g_*$ preserves sheaves in the Nisnevich 
or \'etale topology.     

Since $(g^{-1})_* g_* F = F$, we have 
\begin{equation}\label{adjunction}
\Hom_{PreSh(Cor_k)}\left(g_*F_1,F_2\right) = \Hom_{PreSh(Cor_k)}\left(F_1,(g^{-1})_* F_2\right). 
\end{equation} 

For every smooth scheme $X$ over $k$ the diagonal $\Delta \in X_g \times X_g$ gives a natural morphism 
\eq{natg}
{
\psi_g(X):\Z_{tr}(X) \xr{} (g^{-1})_*\Z_{tr}(X_{g^{-1}}).
}
Here, a section $Z\in Cor(Y,X)$ maps to $(g\times g)_*Z\in Cor(Y_{g^{-1}},X_{g^{-1}})$, and 
natural means functorial in $X$, i.e.~for every morphism $f:X\xr{} Y$ the following 
diagram is commutative
$$
\xymatrix{
\Z_{tr}(X) \ar[r]\ar[d]_{\Z_{tr}(f)} & (g^{-1})_*\Z_{tr}(X_{g^{-1}}) \ar[d]^{(g^{-1})_*
\Z_{tr}(g^{-1}(f))} \\
\Z_{tr}(Y) \ar[r] & (g^{-1})_*\Z_{tr}(Y_{g^{-1}}). }
$$ 
If $X\cong X_0\times_{k_0} k$ then $X_g \cong X$ (over $k$) and we get 
\begin{equation}\label{goperationZX}
\Z_{tr}(X) \xr{\psi_g(X)} (g^{-1})_*\Z_{tr}(X_{g^{-1}}) \xr{\cong} (g^{-1})_*\Z_{tr}(X).
\end{equation}
By abuse of notation we write $\psi_g(X)$ for this map, although it depends on 
the isomorphism $X\cong X_0\times_{k_0} k$. The operation extends in the obvious way to the
Suslin complex $C_*(\Z_{tr}(X))$; explicitly, 
\begin{equation} \label{operationSuslin}
\begin{split}
C_i(\Z_{tr}(X))(Y) &\xr{}  (g^{-1})_*C_i(\Z_{tr}(X))(Y) \\
\Delta^i \times Y \times X \supset Z &\mapsto (g\times g\times g)_*Z.
\end{split}
\end{equation}
The following diagram is commutative 
$$
\xymatrix
{
C_*(\Z_{tr}(X)) \ar[r] & (g^{-1})_*C_*(\Z_{tr}(X)) \\
\Z_{tr}(X) \ar[u] \ar[r]^-{\psi_g(X)} & (g^{-1})_*\Z_{tr}(X). \ar[u]
}
$$

Let the arrangement $\mc{A}$ be defined over $k_0$. For a transverse arrangement $\beta$
(defined over $k$) we denote by $g\beta$ the transverse arrangement $A\mapsto g(\beta(A))$. 
By construction of $I_{\beta}$, we see that 
\begin{equation}\label{Itransforms}
\begin{split}
&\psi_g(V-\cup_{A\in \mc{A}}A) \circ I_{\beta}^{-1} =  \\
&(g^{-1})_*I_{g\beta}^{-1} \circ \bigoplus_{M\in \mc{A}^+} 
\psi_g(\beta(M))/\psi_g(\beta(M)-M)\otimes id_{\Hom(S^{\{M\},\{V\}},\Z)}. 
\end{split}
\end{equation}
Note that we identify $\beta(M)_{g^{-1}} \cong g\beta(M)$, so that
$$
\psi_g(\beta(M)): \Z_{tr}(\beta(M)) \xr{} (g^{-1})_*\Z_{tr}(g\beta(M)), 
$$
and similarly for $\psi_g(\beta(M)-M)$.

%%%%<<<<<<<<<<<<< revision 1

\section{Restriction to a linear subspace}
\label{sec3}

\subsection{}
In this section we study the behavior of the decomposition in Theorem \ref{thm:zerlegung} under restriction 
to a linear subspace $V'$ of $V$. We consider the induced arrangement $\mc{A}'=\{M\cap V'; M\in \mc{A}\}$ 
for $V'$ and do not assume any relation between the transverse arrangements  
$\beta'\colon {\mc{A}'}^{+} \xr{} \mc{T}', \beta \colon \mc{A}^{+} \xr{} \mc{T}$. 

By Theorem \ref{thm:zerlegung} the morphism 
$$\Xi\colon \Z_{tr}(V'-\cup_{A'\in \mc{A}'}A')\xr{} \Z_{tr}(V-\cup_{A\in \mc{A}}A)$$
decomposes into parts $\Xi_{M',M}$, where  $(M',M)\in {\mc{A}'}^{+}\times \mc{A}^{+}$. We give a formula for 
the morphism $\Xi_{M',M}$ in \ref{xi}, but cannot say much about it in general. However, the computation suffices 
to describe the restriction $\Het{*}(V-\cup_{A\in \mc{A}}A,\L) \xr{} \Het{*}(V'-\cup_{A'\in \mc{A}'}A,\L)$ for the
\'etale cohomology in finite coefficients $\L$ (Corollary~\ref{satz}), which is the subject of the next section. 

\subsection{}
Let $\beta \colon \mc{A}^{+} \xr{} \mc{T}$ be a transverse arrangement, i.e. $\beta$ satisfies the properties in
Section \ref{betanotation}. For $A\in \mc{A}^+-\{\emptyset\}$ we denote by 
$$
p_A\colon V \xr{\epsilon_A} \beta(A) \times A \xr{} \beta(A)
$$
the projection. To simplify notation we set $0:=\beta(V)$ and consider $V$ as a vector space with zero point $0$.

\begin{definition} 
For all $A \in \mc{A}^+(m)_{\text{nd}}$ we denote by 
$$\psi(A)\in Cor(V, C_m(\beta(A_0)))$$ 
the graph of $ V \times \Delta^m \xr{} \beta(A_0)$ given by 
\nnal
{
(x , (t_0,\dots,t_m)) \mapsto \sum_{i=0}^m t_i \cdot p_{A_i}(x).
}
\end{definition}

In the following we write   
\small$
C^M:=C_*(\Z_{tr}(\beta(M))/\Z_{tr}(\beta(M)-M))
$ \normalsize
for the Suslin complex.

\begin{lemma}
Let $U:=\{N\in \mc{A}^+; N\cap Y \neq \emptyset \}\subset \mc{A}^+$ for a closed subset $Y\subset V$.
If $\mc{N}\subset U$ is a closed set and $M\in \mc{A}^+$ then  
\eq{phim}
{
\phi^{M,\mc{N}}\colon S^{\{M\},\mc{N}} \xr{} C^{M}(Y), \quad [A] \mapsto \psi(A)_{\mid Y}
}
is a morphism of complexes.
\begin{proof}
We have $\partial_0\psi(A)=0$ for every $A\in \mc{A}^+(m)_{\text{nd}}$ such that \linebreak $A_0=M, A_m\in \mc{N}$, since 
$$\sum_{i=1}^{m} t_i p_{A_i}(y)\in \beta(A_{1}) \subset \beta(A_{0})-A_0.$$

We have $\partial_m\psi(A)=0$ if $A_{m-1}\not\in \mc{N}$, because of  
$$
\partial_m \psi(A)\cap (Y \times \Delta^{m-1} \times (\beta(M)\cap M)) \subset 
(Y\cap A_{m-1}) \times \Delta^{m-1} \times \beta(M),
$$
which is proved in the following Lemma \ref{lemma:psischnitt}. Indeed $\mc{N}$ is closed in $U$, so that
$A_{m-1}\not\in \mc{N}$ and $A_{m-1}\subset A_{m}\in \mc{N}$ imply $A_{m-1}\not\in U$, and $Y\cap A_{m-1}=\emptyset$.
Therefore $\partial_m \psi(A)\subset Y \times \Delta^{m-1} \times (\beta(M) - M)$, which proves the claim. 
\end{proof}
\end{lemma}

\begin{lemma} \label{lemma:psischnitt}
For every $A \in \mc{A}^+(m)_{\text{nd}}$ such that $A_0 = M$, we have 
$$
\psi(A)\cap (Y \times \Delta^{m} \times (\beta(M)\cap M)) \subset 
(Y\cap A_{m}) \times \Delta^{m} \times \beta(M).
$$
\begin{proof}
For $N\in \mc{A}^+-\{ \emptyset \}$ we denote by $x_N$ the intersection point $N\cap \beta(N)$.
Note that if $N,N'\in \mc{A}^+$ such that $\emptyset \not= N'\subset N$ then 
$p_{N}\circ p_{N'}=p_{N}$ and $p_N(x_{N'})=x_N$. Now
$$
\sum_{i=0}^m t_i \cdot p_{A_i}(y)=x_{M}
$$
implies by projection via the map $p_{A_m}$ the identity $p_{A_m}(y)=x_{A_m}$, and therefore $y \in A_m$.
\end{proof}
\end{lemma}

By Lemma \ref{lemma:psischnitt} the morphism $\phi^{M,\mc{N}}$ from \ref{phim} factors through 
\eq{phim2}
{
S^{\{M\},\mc{N}}\otimes \Z_{tr}(Y)/\Z_{tr}(Y- \bigcup_{N\in \mc{N}} N) \xr{} C^{M}.
}

\subsection{}
Let $\mc{M}=\{ M \}, \mc{N}=\mc{A}^+, \mc{P}=\{ V \}$. By \ref{kom} and \ref{phim} we see 
$$
\Hom(S^{\mc{N},\mc{P}},\Z) \xr{} C^M(V) \otimes \Hom(S^{\mc{M},\mc{P}},\Z)
$$
or, equivalently,  
\eq{inv}
{
\Hom(S^{\mc{N},\mc{P}},\Z)\otimes \Z_{tr}(V) \xr{} C^M \otimes \Hom(S^{\mc{M},\mc{P}},\Z). 
}

\begin{proposition} \label{lemma:faktor}
The morphism \ref{inv} factors through 
\eq{}
{
\eta^M \colon S(V,\mc{A}) \xr{} C^M \otimes \Hom(S^{\{M\},\{V\}},\Z) 
}
and the following diagram is commutative
$$
\xymatrix
{
S(V,\mc{A}) \ar[r]^{\eta^M}
&
C^M \otimes \Hom(S^{\{M\},\{V\}},\Z)
\\
S_{\beta}(V,\mc{A}) \ar[u]^{\text{\ref{SbetaS}}} \ar[r]
&
\Z_{tr}(\beta(M))/\Z_{tr}(\beta(M)-M) \otimes  \Hom(S^{\{M\},\{V\}},\Z), \ar[u]_{(\text{\ref{natMorphFCF}})\otimes {\rm id}}
}
$$
where the second horizontal arrow is the projection for the decomposition in \ref{ZerlegungSbeta}. 
\begin{proof}
The morphism \ref{inv} factors through $S(V,\mc{A})$: we need to show that if
$[M \xr{} \dots \xr{} A_{m_1} \xr{} \dots \xr{} A_{m_2}]^*\in \Hom(S^{\{M\},\mc{P}},\Z)$ then
\nnml
{
\psi(M \xr{} \dots \xr{} A_{m_1})\cap (V-A_{m_1})\times \Delta^{m_1}\times \beta(M)\subset \\
V \times \Delta^{m_1}\times (\beta(M)-M),
}
which follows from Lemma \ref{lemma:psischnitt}.

In order to compute the composition with $S_{\beta}(V,\mc{A}) \xr{} S(V,\mc{A})$, we observe that 
\nnml
{
\psi(M \xr{} \dots \xr{} A_{m_1})\cap \beta(A_{m_1})\times \Delta^{m_1}\times \beta(M)\subset \\
\beta(A_{m_1}) \times \Delta^{m_1}\times \beta(A_{m_1}).
}
If $m_1\neq 0$ then $\beta(A_{m_1})\subset \beta(M)-M$ and the restriction to $\beta(A_{m_1})$ vanishes. 
In the case $m_1 = 0$,  $\psi(M)$ induces the identity. This proves the assertion.
\end{proof}     
\end{proposition}

The morphism   
\eq{tauM}
{
\tau^M\colon \Z_{tr}(\beta(M))/\Z_{tr}(\beta(M)-M) \xr{} C^M 
}
is an $\A^1$-weak equivalence \cite[Lemma~14.4]{Voevodskylecture} and Proposition \ref{lemma:faktor} 
implies that $\oplus_{M\in \mc{A}^+}(\tau^M)^{-1} \eta^M$ is the inverse of 
$
S_{\beta}(V,\mc{A}) \xr{} S(V,\mc{A}).
$

\subsection{}

Let $V'\subset V$ be a linear subspace, we set $\mc{A}':=\{A\cap V'; A \in \mc{A}\}$ 
and $\phi\colon \mc{A}^{+} \xr{} {\mc{A}'}^{+}, A\mapsto A\cap V'$. 

Choose transverse arrangements $\beta'\colon {\mc{A}'}^{+} \xr{} \mc{T}'$ and $\beta\colon \mc{A}^{+} \xr{} \mc{T}$.
Let $(M',M)\in {\mc{A}'}^{+} \times {A}^{+}$ and $\mc{N}=\{ N\in \mc{A}^+ ; N\cap V'=M' \}$. 
Using $\phi$ we get  
\eq{SNM'}
{
S^{\mc{N},\{ V \}} \xr{}  S^{\{M'\},\{ V' \}}
}
and we may define ($F_{M'}:=\Z_{tr}(\beta'(M'))/\Z_{tr}(\beta'(M')-M')$)
\ml{xi}
{
\xi_{M',M}\colon F_{M'} \otimes \Hom(S^{\{M'\},\{V'\}},\Z) \xr{{\rm id}\otimes (\text{\ref{SNM'}})} 
F_{M'} \otimes \Hom(S^{\mc{N},\{V\}},\Z) \\ \xr{\text{id}\otimes \text{(\ref{kom})}}
F_{M'} \otimes  S^{\{M\},\mc{N}} \otimes \Hom(S^{\{M\},\{V\}},\Z) \\ \xr{\text{(\ref{phim2})}\otimes \text{id}} 
C^{M}\otimes \Hom(S^{\{M\},\{V\}},\Z);
}
for the last arrow we apply \ref{phim2} to $Y=\beta'(M')$.

\begin{proposition} \label{proposition:transformation}
The morphism $\Z_{tr}(V'-\cup_{A'\in \mc{A}'}A')\xr{} \Z_{tr}(V-\cup_{A\in \mc{A}}A)$ 
in $\mathbf{D}^-[W_{\A}^{-1}]$ is via the decomposition from Theorem \ref{thm:zerlegung} given by 
$$\bigoplus_{(M',M)\in {\mc{A}'}^{+}\times \mc{A}^{+}} (\tau^{M})^{-1}\xi_{M',M},$$ 
where $\tau^M$ is as in \ref{tauM} and $\xi_{M',M}$ is defined in \ref{xi}.
\begin{proof}
Let $f\in Cor(V',V)$ be the inclusion, then we have 
$$
S(f,\phi)\colon S(V',\mc{A}') \xr{} S(V,\mc{A})
$$
from \ref{funkR}. By \ref{nat} it is enough to show that $\xi_{M',M}$ is given by 
\nnml %{M'M}
{
F_{M'} \otimes  \Hom(S^{\{M'\},\{V'\}},\Z) \xr{} S_{\beta'}(V',\mc{A}') 
\xr{\text{\ref{SbetaS}}} S(V',\mc{A}') \xr{} S(V,\mc{A}) \\
\xr{\text{Proposition \ref{lemma:faktor}}} C^M \otimes \Hom(S^{\{M\},\{V\}},\Z).
} 
This verification can be carried out directly.
\end{proof}
\end{proposition}

\section{\'Etale cohomology}
\label{sec4}

\subsection{}
Let $V$ and $\mc{A}$ be as in Section \ref{betanotation}, we set $X:=V-\cup_{A\in \mc{A}} A$. 
In order to compute the \'etale cohomology ring  
$H^*_{\text{\'et}}(X,\L)$, where \linebreak $\L=\Z/n$ and $n$ is invertible in the ground field $k$, we consider 
the derived category $\mathbf{D}^-=\mathbf{D}^-(Sh_{\text{\'et}}(Cor_k,\Z/n))$ of \'etale
sheaves of $\Z/n$-modules with transfers and the localization $\DMet=\mathbf{D}^-[W^{-1}_{\A}]$. We use 
the natural isomorphism
\eq{DMetD}
{
\Hom_{\DMet}(\Ltr(U),\L[p])\cong \Hom_{\mathbf{D}^-}(\Ltr(U),\L[p]) \cong \Het{p}(U,\L)
}  
for all $k$-schemes $U$ and integers $p$ \cite[Proposition~10.7]{Voevodskylecture}. 
Furthermore we have $\Lambda(q)\cong \mu_{n}^{\otimes q}$ in $\DMet$, where $\mu_{n}$ is the locally constant sheaf of $n$-th 
roots of unity \cite[Proposition~10.6]{Voevodskylecture}.

%%%>>>>>>>>>> revision 1
\subsection{}
The isomorphism $I_{\beta}$ from Theorem \ref{thm:zerlegung} takes for finite coefficients $\Lambda$ 
and the \'etale topology the form 
\begin{equation} \label{Lambdazerlegung}
\Lambda_{tr}(X) \xr{I_{\beta}}  
 \bigoplus_{M\in \mc{A}^+-\{\emptyset\}} \mu_n^{\otimes \cd M}[2\cd M] \otimes  \Hom(S^{\{M\},\{V\}}\otimes \Lambda,\Lambda) 
\end{equation}
where $\cd M:=\codim M$.  Although the decomposition depends  a priori on the transverse 
arrangement, we denote by $h_M(V,\mc{A})$ the direct summand for $M$.

\begin{thm} \label{satzhomological} 
Let $k$ be an algebraically closed field and $V'\subset V$ a linear subspace,  
$\mc{A}'=\{A\cap V'; A\in \mc{A}\}, X'=V'-\cup_{A'\in \mc{A}'} A',$ and 
$\phi\colon \mc{A}^{+} \xr{} {\mc{A}'}^{+}, A\mapsto A\cap V'$. 
\begin{enumerate}
\item[(i)] The morphism \ref{Lambdazerlegung} 
does not depend on the choice of a transverse arrangement.  
\item[(ii)] The composite   
$$
h_{M'}(V',\mc{A}') \xr{} \Lambda_{tr}(X') \xr{} \Lambda_{tr}(X) \xr{} h_{M}(V,\mc{A})
$$
vanishes, if $M\cap V'\neq M'$ or $\cd_{M'}V'\neq \cd_{M} V$. Otherwise, it is induced by
the dual of 
$$
S^{\{M\},\{V\}} \xr{} S^{\{M\cap V'\},\{V'\}}; \quad [A] \mapsto [\phi(A)].
$$
\item[(iii)] If the arrangement is defined over $k_0\subset k$ and $g\in \Aut_{k_0}(k)$,  
then the morphism
$$
h_{M}(V,\mc{A}) \xr{} \Lambda_{tr}(X) \xr{} (g^{-1})_*\Lambda_{tr}(X) \xr{} (g^{-1})_*h_{M'}(V,\mc{A}) \xr{} h_{M'}(V,\mc{A}) 
$$
vanishes if $M\neq M'$, and otherwise it is $g\otimes id_{\Hom(S^{\{M\},\{V\}},\Lambda)}$ with 
$g:\mu_{n}^{\otimes \cd M}\xr{} \mu_{n}^{\otimes \cd M}$ 
coming from the action of $\Aut_{k_0}(k)$ on $\mu_n(k)$.  
\end{enumerate}
\begin{proof}
We choose the  transverse arrangements $\beta$ in $V$ and $\beta'$ in $V'$ independently, so that (i) follows from (ii) for $V'=V$. 

(ii) By Proposition \ref{proposition:transformation} and \ref{DMetD} we have to consider $\xi_{M',M}$ in $\DMet$. 
Obviously we have $\xi_{M',M}=0$ if $M'=\emptyset$ or $M=\emptyset$. 

For $\mc{N}=\{ N\in \mc{A}^+ ; N\cap V'=M' \}$ the morphism 
\small
$$
S^{\{M\},\mc{N}}\otimes \Z_{tr}(\beta'(M'))/\Z_{tr}(\beta'(M')-M') \xr{}  C_*(\Z_{tr}(\beta(M))/\Z_{tr}(\beta(M)-M))
$$\normalsize
from \ref{phim2} corresponds to a morphism of complexes of $\L$-modules 
\eq{phim2khomological}
{
S^{\{M\},\mc{N}}\otimes \mu_{n}^{\otimes \cd{M'}}(k)[2\cd{M'}] \xr{} \mu_{n}^{\otimes \cd{M}}(k)[2\cd{M}].
}  
The complex $S^{\{M\},\mc{N}}$ is concentrated in the degrees $-\cd M+c,\dots,0$, where $c:={\rm min}\{ \codim N ; N \in \mc{N}\}$,
and we have $c\geq \cd{M'}$. 
Therefore \ref{phim2khomological} vanishes if $\cd{M'}\neq \cd{M}$. 

In the case $ \cd{M'} = \cd M $ we see $M'=M\cap V'$ and  $S^{\{M\},\mc{N}}=\L$, 
thus \ref{phim2khomological} is the identity map. 
Now it is easy to verify the assertion of (ii) by using the definition \ref{xi} of $\xi_{M',M}$. 

(iii) By using (i) the vanishing for $M\neq M'$ follows immediately from formula \ref{Itransforms}.
Formula \ref{Itransforms} also implies that the general case can be reduced to 
$\mc{A}=\{M,V\}$ with $M$ of codimension $c=\cd M$ ($M$ defined over $k_0$);
by homotopy invariance we may assume that $M=\{0\}$ and $V=\A^{c}$. 
We need to compute the composite 
\nnml
{
\mu_{n}^{\otimes c}[2c] \xr{} \Lambda_{tr}(V)/\Lambda_{tr}(V-0) \xr{} g_*(\Lambda_{tr}(V)/\Lambda_{tr}(V-0))\xr{} \\ 
g_*\mu_{n}^{\otimes c}[2c]=\mu_{n}^{\otimes c}[2c].
}  
The arguments of Section \ref{Gm} show that we may replace \linebreak 
$\Lambda_{tr}(V)/\Lambda_{tr}(V-0)$
by $\L_{tr}(\G_m^c)[c]$. The morphism $\mu_n^{\otimes c}[c]\xr{} \L_{tr}(\G_m^c)$ is defined by
$$
\mu_n^{\otimes c}[c] \xr{\tau} C_*(\L_{tr}(\G_m^c)) \xl{} \L_{tr}(\G_m^c)
$$  
where $\tau$ is a morphism in $\mathbf{D}^-$; $\tau$ is uniquely determined by its
stalk on $k$ (see \ref{adjoint2} below). For $c=1$ the map $\mu_n(k)\xr{} Cor(\Delta,\G_m)$
is defined by $\zeta\mapsto Z(\zeta)$ with 
\begin{equation*}
%\begin{split}
Z(\zeta)_{(0,1)\times {\G_m}}=n\cdot [\zeta], \quad 
Z(\zeta)_{(1,0)\times {\G_m}}=n\cdot [1].
%\end{split}
\end{equation*}
In coordinates $(t_0,t_1),t_0+t_1=1,$ for $\Delta$,
and $x$ for $\G_m$, we may take $Z(\zeta)=\{x^n+t_0f_0(x)+t_1f_1(x)=0\}$ with 
$f_0=(x-1)^{n}-x^n$ and $f_1=(x-\zeta)^{n}-x^n$.

For $c\geq 1$ the map $\mu_n(k)^{\otimes c}\xr{} Cor(\Delta^c,\G_m^c)$ is 
$$\zeta_{1}\otimes \dots \otimes \zeta_c\mapsto Z(\zeta_1)\times \dots \times Z(\zeta_c),$$
and therefore the composite 
$$\mu_n(k)^{\otimes c}[c]\xr{} C_*(\L_{tr}(\G_m^c))(k) \xr{} (g^{-1})_*C_*(\L_{tr}(\G_m^c))(k)$$
maps $\zeta_{1}\otimes \dots \otimes \zeta_c$ to $Z(g(\zeta_1))\times \dots \times Z(g(\zeta_c))$ 
(by \ref{operationSuslin}), as claimed. 
\end{proof}
\end{thm}

\subsection{} \label{combinatoric}
We define the decomposition of \'etale cohomology by 
\begin{equation*}
\begin{split}
&\Het{p}(X,\L) \xr{\ref{DMetD}} \Hom_{\DMet}(\Ltr(X),\L[p]) \xr{\ref{Lambdazerlegung}} \\
\bigoplus_{M \in \mc{A}^+} & \Hom_{\DMet}(h_M(V,\mc{A}),\L[p]) \cong 
\bigoplus_{M \in \mc{A}^+} \Hom_{\mathbf{D}^-}(h_M(V,\mc{A}),\L[p]). 
\end{split}
\end{equation*}
The last isomorphism follows from the fact that $L$ is $\A^1$-local \cite[Definition~9.17, Corollary~9.25]{Voevodskylecture}. In order to give an explicit description, we note that if 
$C$ is a complex of constant sheaves and $k$ is algebraically closed, then 
\begin{equation} \label{adjoint2}
\Hom_{\mathbf{D}^-}(C,F)=\Hom_{D^{-}(\Lambda)}(C(k),F(k))
\end{equation}
for every $F\in \mathbf{D}^-$. Here, $D^{-}(\Lambda)$ is the (bounded above) derived category
of $\Lambda$-modules, and $F(k)$ is the complex of $k$-sections (or the stalk of $F$ at $k$). 
Therefore, we get
$$
\Hom^p_{\mathbf{D}^-}(h_M(V,\mc{A}),\L)=H^{p-2\cd M}\left(S^{\{M\},\{V\}}\otimes \L \right)\otimes \Hom(\mu_{n}^{\otimes \cd M}(k),\Lambda).
$$
We denote this direct summand in $\Het{p}(X,\L)$ by $h^p_M(V,\mc{A})$.

\subsection{} \label{GoperationviaDM}
If the arrangement is defined over a subfield $k_0\subset k$ then the 
Galoisgroup $G=\Aut_{k_0}(k)$ acts on the \'etale cohomology of $X$ by
$$
\Het{p}(X,\Lambda)\xr{=} \Het{p}(X,g_*\Lambda) \xr{\cong} \Het{p}(X,\Lambda), 
$$ 
where the first equality follows from $F(X)=(g_*F)(X)$ for every sheaf (and is
true for every sheaf, i.e.~$\Het{p}(X,F)\xr{=} \Het{p}(X,g_*F)$), and the second
isomorphism is induced by the isomorphism $\Lambda\xr{} g_*\Lambda$ of sheaves.

In $\mathbf{D}^-$ the action is given by  
\begin{equation*}
\begin{split}
&\Hom^p_{\mathbf{D}^-}(\Z_{tr}(X),\L) \xr{\ref{goperationZX}} 
\Hom^p_{\mathbf{D}^-}((g^{-1})_*\Z_{tr}(X),\L) \xr{\ref{adjunction}} \\  
&\Hom^p_{\mathbf{D}^-}(\Z_{tr}(X),g_*\L) \xr{} \Hom^p_{\mathbf{D}^-}(\Z_{tr}(X),\L).  
\end{split}
\end{equation*}
In order to prove this, it is sufficient to show that the composite of the first two 
arrows is $\Het{p}(X,\Lambda)\xr{=} \Het{p}(X,g_*\Lambda)$. By using an injective 
resolution for $\Lambda$ this follows immediately from the definitions.

\begin{korollar} \label{satz}  
Let $k$ be an algebraically closed field and $V'\subset V$ a linear subspace,  
$\mc{A}'=\{A\cap V'; A\in \mc{A}\}$ and $\phi\colon \mc{A}^{+} \xr{} {\mc{A}'}^{+}, A\mapsto A\cap V'$. 
Moreover, let $k_0\subset k$ be a subfield such that the arrangement is defined over $k_0$.
We denote by $X=V-\cup_{A\in \mc{A}}A$ (resp. $X'=V'-\cup_{A'\in \mc{A}}A'$) the complement
of the arrangement.
\begin{enumerate}
\item[(i)] The decomposition 
\begin{equation*}
\begin{split}
\Het{p}(X,\L)&\cong \bigoplus_{M\in \mc{A}^+} h^p_M(V,\mc{A}) \\
\cong \bigoplus_{M\in \mc{A}^+} & H^{p-2\cd M}\left(S^{\{M\},\{V\}}\otimes \L \right)\otimes \Hom(\mu_{n}^{\otimes \cd M}(k),\Lambda).
\end{split}
\end{equation*}
does not depend on the choice of a transverse arrangement. It is a decomposition of 
$G=Gal(k/k_0)$-modules, where the $G$-action on 
$H^{p-2\cd M}\left(S^{\{M\},\{V\}}\otimes \L \right)$ is 
trivial.
\item[(ii)] For the pullback  
$
\iota^* \colon \Het{*}(X,\L) \xr{} \Het{*}(X',\L)
$
we have 
$$\iota^*(h^*_M(V,\mc{A}))=0$$
if $M\cap V'=\emptyset$ or $\codim_{V'}(M\cap V')<\codim_{V}(M)$. 
In the case $\codim(M\cap V')=\codim(M)$ we have 
$$
\iota^*(h^*_M(V,\mc{A})) \subset h^*_{M\cap V'}(V',\mc{A}')
$$
and the restriction of $\iota^*$ to $h^*_M(V,\mc{A})$ is induced by 
$$
S^{\{M\},\{V\}} \xr{} S^{\{M\cap V'\},\{V'\}}, \quad [A] \mapsto [\phi(A)].
$$ 
\end{enumerate}
\begin{proof}
(i) By Theorem \ref{satzhomological}(i) the decomposition doesn't depend on the choice 
of a transverse arrangement. The combinatorial description has been proved in 
Section \ref{combinatoric}. By using Section \ref{GoperationviaDM} the Galois operation follows
from Theorem \ref{satzhomological}(iii) by a straightforward calculation. \\ 
(ii) Follows from Theorem \ref{satzhomological}(ii).
\end{proof}
\end{korollar}

%%%%<<<<<<<<< revision 1

\subsection{External cup product}
Let $V,Y$ be smooth $k$-schemes and $\mc{A}$ (resp. $\mc{B}$) a set of closed subsets  
of $V$ (resp. $Y$), as in Section \ref{notation}. If 
$\mc{M}_1\subset \mc{A}^+$ and $\mc{M}_2\subset \mc{B}^+$ are locally closed subsets then there is a morphism 
\eq{extcup}
{
\Hom(S^{\mc{M}_1,\{V\}},\L) \otimes \Hom(S^{\mc{M}_2,\{Y\}},\L) \xr{} \Hom(S^{\mc{M}_1\times \mc{M}_2,\{V\times Y\}},\L)
}
defined by 
\small \begin{multline} \label{extcupformula}
[A_0 \xr{} \dots \xr{} A_{m-1}\xr{} V]^* \otimes [B_0\xr{} \dots \xr{} B_{n-1}\xr{} Y]^* 
\mapsto \\
\sum
[A_0\times C_0 \xr{} A_1\times C_1 \xr{} \dots \xr{} V\times B_0 \xr{} V\times B_1 \xr{} \dots \xr{} V\times Y]^* \end{multline} \normalsize
where the sum is over all $C_0\xr{} \dots \xr{} C_{m-1}\xr{} B_0 \in \mc{B}^+(m)$ such that $C_0\in \mc{M}_2$ (and therefore $C_i\in \mc{M}_2$ for every $i$). 
The morphism \ref{extcup} can be interpreted as follows. 
Write $\mc{M}_1=U_1\cap Z_1$ resp. $\mc{M}_2=U_2\cap Z_2$ with open sets $U_1,U_2$ and closed sets 
$Z_1,Z_2$. From the obvious functors $U_1\times U_2\xr{p_1} U_1$ resp. $U_1\times U_2\xr{p_2} U_2$ we have  
\eq{Bproduct}
{
\abs{U_1 \times U_2} \xr{(p_1,p_2)} \abs{U_1}\times \abs{U_2}
}
and \ref{extcup} gives $p^*_1\cup p^*_2$ on the relative cohomology:   
\nnml
{
H^*(\abs{U_1},\abs{U_1-\{V\}}\cup \abs{U_1-\mc{M}_1}) \otimes H^*(\abs{U_2},\abs{U_2-\{Y\}}\cup \abs{U_2-\mc{M}_2}) \xr{} \\
H^*(\abs{U_1 \times U_2}, \abs{U_1\times U_2 - \{ V\times Y \}} \cup \abs{U_1\times U_2 - \mc{M}_1\times \mc{M}_2}).
}
By the Eilenberg-Zilber Theorem we know 
$$
H^*(\abs{\mc{D}_1 \times \mc{D}_2};\L) \cong H^*(\abs{\mc{D}_1} \times \abs{\mc{D}_2};\L).
$$
Using the Kuenneth-isomorphism and the short exact sequences from \ref{SMNaufloesung} we conclude that \ref{extcup}
is a quasi-isomorphism.  

If $\mc{M}_1=\mc{A}^+, \mc{M}_2=\mc{B}^+$ we use \ref{extcup} to define the right vertical arrow in the commutative 
diagram  
\eq{dia}
{
\xymatrix
{
\Ltr(V-\cup_{A\in \mc{A}}A) \otimes^{tr} \Ltr(Y-\cup_{B\in \mc{B}}B) \ar[r]^-{\text{quis}} \ar[d]^{\cong}
&
S(V,\mc{A})\otimes^{tr}_{\mathbb{L}} S(Y,\mc{B}) \ar[d]
\\
\Ltr(V\times Y -\cup_{(A,B)\in \mc{A}\times \mc{B}} A\times B)  \ar[r]^-{\text{quis}}
&
S(V\times Y,\mc{A}\times \mc{B}).
}
}

\subsection{}
Via the isomorphism \ref{DMetD} the external product on the  \'etale cohomology for two smooth $k$-schemes  
$U_1,U_2$ corresponds to  
\small $$ 
\xymatrix @C-3pc @R-1pc
{
\Hom_{\DMet}^{p_1}(\Ltr(U_1),\L) \otimes \Hom^{p_2}_{\DMet}(\Ltr(U_2),\L) \ar[d]
\\ 
\Hom^{p_1+p_2}_{\DMet}(\Ltr(U_1)\otimes^{tr} \Ltr(U_2),\L) \ar[r]^-{=} & \Hom^{p_1+p_2}_{\DMet}(\Ltr(U_1 \times U_2),\L).
} $$ \normalsize
From the diagram \ref{dia} it is easy to see that
\small$$
\Het{*}(V-\cup_{A\in \mc{A}}A,\L)\otimes \Het{*}(V-\cup_{A\in \mc{A}}A,\L) \xr{} 
\Het{*}(V\times V -\cup_{(A,B)\in \mc{A}\times \mc{A}}A\times B,\L)
$$\normalsize 
preserves the decomposition of Corollary \ref{satz}, i.e. 
$$
h^*_{M_1}(V,\mc{A})\otimes h^*_{M_2}(V,\mc{A}) \xr{}  h^*_{M_1\times M_2}(V\times V,\mc{A}\times \mc{A})
$$ 
for every $M_1,M_2\in \mc{A}^{+}-\{\emptyset\}$. Explicitly the map is given by  
\small\ml{extcupet}
{ 
H^{*}(S^{\{M_1\},\{V\}}\otimes \L[-2\cd{M_1}])
\otimes H^{*}(S^{\{M_2\},\{V\}}\otimes \L [-2\cd{M_2}]) \\
\xr{} H^{*}(S^{\{M_1\},\{V\}}\otimes S^{\{M_2\},\{V\}} \otimes  \L[-2(\cd{M_1}+\cd{M_1})]) \\
\xr{\cong} H^{*}(S^{\{M_1\times M_2\},\{V\times V\}} \otimes  \L[-2(\cd{M_1}+\cd{M_1})])
}\normalsize
and twist with $\mu_n^{\otimes -\cd M_1 - \cd M_2}$. Here the last map comes from the quasi-isomorphism \ref{extcup}. 

\subsection{} 
The multiplication on the cohomology ring is the composition of the external product and the pullback to the diagonal. 
From the computation of the external product \ref{extcup} we get the following corollary of Corollary \ref{satz}.

\begin{korollar} \label{Multiplikation}
The  cup product of $\Het{*}(V-\cup_{A\in \mc{A}}A,\L)$ vanishes on \\  
$h^*_{M_1}(V,\mc{A})\otimes h^*_{M_2}(V,\mc{A})$ 
if $M_1\cap M_2=\emptyset$ or $\cd M_1+\cd M_2 > \cd (M_1\cap M_2)$. 
If $\cd M_1+\cd M_2 =\cd (M_1\cap M_2)$ then $h^*_{M_1}(V,\mc{A})\otimes h^*_{M_2}(V,\mc{A})$  maps to $h^*_{M_1\cap M_2}(V,\mc{A})$. 
Explicitly the map can be computed as the composition of the external product \ref{extcupet} and the map   
\nneq{
H^{*}(S^{ \{M_1\times M_2\},\{V\times V\} } \otimes  \mu_n(-c)[-2c]) \xr{}  
H^{*}(S^{\{M_1\cap M_2\},\{V\}} \otimes  \mu_n(-c)[-2c])
}
($c:=\cd M_1+\cd M_2$) induced by  
$\mc{A}\times \mc{A}\xr{} \mc{A}, A_1\times A_2 \mapsto A_1\cap A_2$.
\end{korollar} 

\section{Comparison with known decompositions}

\subsection{}\label{defprop}
%%%>>>>>>>>>>>< revision2
 In this section we work with $\Lambda=\Z/n$ and $n^{-1}\in k$.
%%%<<<<<<<<<<<< revision2
Assume that for all affine spaces $V$ and finite sets $\mc{A}$ of linear subspaces  
(such that $A_1,A_2\in \mc{A}\Rightarrow A_1\cap A_2 \in \mc{A}$) we have a decomposition  
\eq{deco}
{
\Het{*}(V-\cup_{A\in \mc{A}}A,\L) = \bigoplus_{A\in \mc{A}^{+} } h_A(V,\mc{A}).
}
We assume $h_{\emptyset}(V,\mc{A})=0$ and consider the following properties:
\begin{enumerate}
\item[F1:] For every affine-linear morphism $\iota\colon V'\xr{} V$ we have   
$$
\iota^{*} h_A(V,\mc{A}) \subset h_{\iota^{-1}(A)}(V',\mc{A}'), 
$$ 
where $\mc{A}':=\{ \iota^{-1}(A) ; A\in \mc{A}\}$.
\item[F2:] For all $A'\in \mc{A}$ and $\mc{A}' \subset \mc{A}$ defined by \small$\mc{A}' := \{A \in \mc{A}; A \supset A' \}$ \normalsize 
the map  $ j^*\colon \Het{*}(V-\cup_{A\in \mc{A}'}  A,\L)  \xr{} \Het{*}(V-\cup_{A\in \mc{A}} A,\L)$ 
preserves the decomposition and furthermore 
 $$
 j^* h_{A''}(V,\mc{A}')= h_{A''}(V,\mc{A})
 $$  
for every $A''\in \mc{A}'$.
\end{enumerate}  

\begin{proposition} \label{propeind}
\begin{enumerate}
\item[(i)] 
If $h,h'$ are two decompositions satisfying (F2), then 
$$
h'_A(V,\mc{A}) \xr{} \Het{*}(V-\cup_{A\in \mc{A}}A,\L) \xr{} h_A(V,\mc{A})
$$ 
is an isomorphism.
\item[(ii)] Properties (F1),(F2) determine the decomposition uniquely.
\end{enumerate}
\begin{proof}
\emph{Proof of (i).} By property (F2) it is enough to consider $\mc{A}'= \{A'\in \mc{A}, A'\supset A\}$ and 
$V-\cup_{A'\in \mc{A}'} A'$. For every $A'\in \mc{A}'$ we have 
$h'_{A'}(V,\mc{A}')\subset \oplus_{A''\supset A'} h_{A''}(V,\mc{A}')$; this implies the assertion.

\emph{Proof of (ii).} We proceed by induction on the dimension of $V$. By (F2) we may assume that $A$ is the minimal 
element of the arrangement $\mc{A}$, and we need to prove that $h'_A(V,\mc{A})=h_A(V,\mc{A})$ as subgroups of  
$\Het{*}(X,\L)$ (as usual, $X:=V-\cup_{A'\in \mc{A}}A'$).   

For every $A'\in \mc{A}-\{A\}$ we choose a linear subspace $V'$ transverse to $A'$, and such that $A''\supset V'\cap A'$ implies 
$A''\supset A'$ for all $A''\in \mc{A}$.
We set $\mc{A}':= \{A''\in \mc{A}; A''\supset A'\}$ and $U:=V-\cup_{A''\in \mc{A}'}A''$, and consider 
the commutative diagram
$$
\xymatrix
{
V'\cap U \ar[r]
&
U 
\\
V'\cap X  \ar[u] \ar[r]^-{\imath}
&
X. \ar[u]_{\jmath} 
}
$$
In the composition 
\nneq
{
\tau\colon \Het{*}(X,\L) \xr{\imath^{*}} \Het{*}(V'\cap X,\L) \xr{}  
\Het{*}(V'\cap U,\L) \xr{(pr_{A'}^*)^{-1}}  \Het{*}(U,\L), 
}
the second arrow depends on the decomposition of the cohomology of  $V'\cap X$, 
which is unique by induction. Here, we denote by $pr_{A'}\colon V \xr{} V'$ the projection along $A'$. 
It is easy to see that 
\small$$
\Het{*}(X,\L) \xr{\tau} \Het{*}(U,\L) \xr{\jmath} 
\Het{*}(X,\L)
$$ \normalsize
is the projector
%% [Es wird hier benutzt, dass 
%% $$ 
%% V-\cup_{A''\in \mc{A}''} A'' \xr{pr_{A'}} V'-\cup_{A''\in \mc{A}''}(A''\cap V') \xr{} V-\cup_{A''\in \mc{A}''} A''
%% $$ 
%% die Identitaet auf der Kohomologie induziert, was wegen 
%% $$(V'-\cup_{A''} (A''\cap V'))\times A'=V'-\cup_{A''}(A''\cap V')$$
%% klar ist]
of the direct summand $\oplus_{A''\supset A'} h_{A''}(V,\mc{A})$. By the same argument it is the projector of 
$\oplus_{A''\supset A'} h'_{A''}(V,\mc{A})$. Using the equality of these projectors for every $A'\supsetneq A$ we see
that the projectors of $h_{A'}(V,\mc{A})$ and $h'_{A'}(V,\mc{A})$ are equal for all  $A'\in \mc{A}-\{A\}$; we denote the projector by $P_{A'}$. 
Since ${\rm id} - \sum_{A'\supsetneq A} P_{A'}$ is the projector for $A$, we are done. 
\end{proof}
\end{proposition}

\subsection{}
We need to show that the decomposition of Corollary~\ref{satz}(i) satisfies (F1),(F2). 
\emph{Proof of (F1)}: A linear morphism $\iota \colon V'\xr{} V$ may be written as the composition of a projection and an inclusion.
For the inclusion the claim follows from Corollary~\ref{satz}(ii). If $\iota \colon V'\xr{} V$ is a projection, we may use a section
$s\colon V\xr{} V'$ to define a transverse arrangement $\beta'$ for $\mc{A}'$:
$$
\beta'(\iota^{-1}(A)):=s(\beta(A)) \quad \text{for all $A\in \mc{A}$.}
$$ 
The arrangement is compatible with $\beta$ in the sense of \ref{seckompbb'}. The assertion follows from the  
commutative diagram \ref{diakompbb'}. More precisely, $\iota$ preserves the motivic decomposition  
of Theorem~\ref{thm:zerlegung} for $\beta'$ and $\beta$. \\
\emph{Proof of (F2)}: Since $\mc{A}'\subset \mc{A}$ every transverse arrangement $\beta$ for $\mc{A}$ is one for $\mc{A}'$ 
too. One can argue as in the prove of (F1).

\subsection{}
In \cite[Prop.~2.2]{DGM} it is proven that the decomposition of \cite{DGM} satisfies (F2). So the direct summands 
of the two decompositions are isomorphic by Proposition~\ref{propeind}(i). If the decomposition of \cite{DGM} satisfies
\cite[Lemma~3.3]{DGM} then  \cite[Theorem~3.2]{DGM} holds for \'etale cohomology too. This implies (F1), and the two 
decompositions would agree by Proposition~\ref{propeind}(ii).  

%%%%>>>>>>>>>>>>>> revision 2

\subsection{}
In \cite[Exemple~1.14]{DGM} the decomposition of the 
compact \'etale cohomology of $X$ is given by 
\begin{equation} \label{DGM1.14}
\bigoplus_{A\in \mc{A}^{+}} H^{*-2\dim A}_{\{A\}}([A,V],\epsilon_{V!} \Lambda) \otimes \mu_n(-\dim A) \xr{\cong} H^*_c(X,\Lambda), 
\end{equation}
where $[A,V]$ is equipped with the topology from \ref{arrangementopology}, and 
$\epsilon_V:\{V\}\xr{} [A,V]$ is the inclusion of the point $V$. In order to clarify 
the connection with our decomposition we need to compute the cohomology groups of 
\ref{DGM1.14} more explicitly.

\subsection{} Let $\mc{C}$ be a category with finitely many objects and morphism, 
with the topology from \ref{arrangementopology} 
(only $\mc{C}=[A,V]\subset \mc{A}^+$ will be important for us). 
For a sheaf $F$ of $\Lambda$-modules 
on $\mc{C}$ we denote by $F_x$ the stalk 
at $x\in \mc{C}$. We have $F_x=F(U(x))$ with  
$U(x)=\{ y\in \mc{C}; \exists x\xr{} y \}$ the smallest open set that contains $x$. 
Denote by $\epsilon_x$ the inclusion of a point $x$ in $\mc{C}$.  
We define a complex of sheaves $G(F,\mc{C})$ on $\mc{C}$ by setting 
$$
G(F,\mc{C})^m = \bigoplus_{x_0\xr{} \dots \xr{} x_m \in \mc{C}(m)_{{\rm nd}}} 
\epsilon_{x_0*} F_{x_m}.
$$
Note that $\epsilon_{x_0*} F_{x_m}$ is the extension by zero of the constant sheaf with fibre 
$F_{x_m}$ on the closure of $x_0$. The differential is defined by $\sum_i (-1)^i \partial^i$
with 
$$
{\rm pr}_x \circ \partial^i = \begin{cases} (\epsilon_{x_1*}F_{x_{m}} \xr{}  \epsilon_{x_0*}F_{x_{m}}) \circ {\rm pr}_{\partial_i(x)} & \text{if $i=0$,}  \\
\epsilon_{x_0*}(F_{x_{m-1}}\xr{} F_{x_m}) \circ  {\rm pr}_{\partial_i(x)}  & \text{if $i=m$,}  \\
\epsilon_{x_0*}(id_{F_{x_{m}}})  \circ  {\rm pr}_{\partial_i(x)} & \text{otherwise,}  \\
\end{cases}
$$ 
for all $x\in  \mc{C}(m)_{{\rm nd}}$. For $i=0$ the morphism is induced by 
$\epsilon_{x_0}^*\epsilon_{x_1*}(F_{x_{m}}) \xr{=} F_{x_m}$, and for $i=m$ by
the restriction 
$$
F_{x_{m-1}} = F_{U(x_{m-1})} \xr{}  F_{U(x_{m})} = F_{x_m}, 
$$
because the arrow $x_{m-1}\xr{} x_m$ implies $U(x_{m}) \subset U(x_{m-1})$.

Since $G(F,\mc{C})^0=\oplus_{x\in \mc{C}} \epsilon_{x_*}F_x$ there is a natural
morphism $F\xr{} G(F,\mc{C})$ which maps a section $s\in F(U)$ to $\oplus_{x\in U} s_x$.

\begin{lemma} \label{lemmaFquisG}
The map
$$
F\xr{} G(F,\mc{C})
$$
is an quasi-isomorphism.
\begin{proof}
The complex $G(.,\mc{C})$ is functorial and preserves exact sequences. By the D\'evissage
Lemma \cite[1.3]{DGM} it suffices to prove the statement for $F=\epsilon_{x*}M$ with $x\in \mc{C}$
and a $\Lambda$-module $M$.

For $y\in \mc{C}$ the components $G(F,\mc{C})^m_y$ are given by 
$$
\bigoplus_{x_0\xr{} \dots \xr{} x_m} M 
$$
where the sum is over  $x_0\xr{} \dots \xr{} x_m \in \mc{C}(m)_{{\rm nd}}$ such that
$\exists y \xr{} x_0$ and $\exists x_m\xr{} x$ holds. 
Therefore $G(F,\mc{C})_y$ computes
the cohomology of the classifying space of the full subcategory 
$\mc{C}_{y\xr{}x}=\{z;\exists y \xr{} z\xr{} x\}$ of $\mc{C}$. 
This space is empty if there is no map $y\xr{} x$, in this case we see 
$G(F,\mc{C})_y=0$ and $F_y=0$, so that the claim holds. 

Now assume that there is a map $y\xr{} x$, then $y$ is contained in the closure of $x$,  
thus $F_y=M$. The classifying space of
$\mc{C}_{y\xr{}x}$ is contractible because it has a minimal element $y$, and
therefore $$H^*(G(F,\mc{C})_y)=H^0(G(F,\mc{C})_y)=M.$$  
\end{proof}
\end{lemma}

We work with the resolution $G(F,\mc{C})$ because the sheaves $\epsilon_{x*} M$ are 
acyclic for the functors we shall need.

\begin{lemma} \label{resolutionG}
For any $\Lambda$-module $M$ and any point $x\in \mc{C}$ the sheaf $F=\epsilon_{x*}M$ 
is for the following functors acyclic:
\begin{itemize}
\item[(i)] $\pi_*$ with $\pi:\mc{C} \xr{} \mc{C}'$ a functor,
\item[(ii)] $\Gamma_{\mc{C}'}$, i.e. global sections with support in $\mc{C}'$, 
where $\mc{C}'\subset \mc{C}$ is a closed set.
\end{itemize}
\begin{proof} \label{acyclic}
In order to prove (i) we need to show that 
$$
R^{*}\pi(F)_y = H^i(\pi^{-1}U(y), F_{\mid \pi^{-1}U(y)})
$$
vanishes if $i>0$. If $x\in \pi^{-1}U(y)$ then we denote by $\epsilon_{x}'$ the inclusion 
of $x$ in $\pi^{-1}U(y)$.  We have
$$ 
\epsilon_{x*}M_{\mid \pi^{-1}U(y)} = \begin{cases} 0 & \text{if $x\not\in \pi^{-1}U(y)$} \\
                                                   \epsilon_{x*}'M & \text{if $x\in\pi^{-1}U(y)$.}
\end{cases}
$$
Since ${\epsilon'_{x}}^*$ and $\epsilon_{x*}'$ are adjoint, $\epsilon_{x*}'M$ is injective if 
$M$ is injective. The functor $\epsilon_{x*}'$ is exact, and therefore
$$
H^i(\pi^{-1}U(y), \epsilon_{x*}M_{\mid \pi^{-1}U(y)}) = H^i(\{x\},M) = 0 \quad \text{if $i>0$.}
$$

For (ii) we use the long exact sequence 
$$
H^*_{\mc{C}'}(\mc{C},F) \xr{} H^*(\mc{C},F) \xr{} H^*(\mc{C}-\mc{C}',F) \xr{+1}
$$
and observe that by (i) we have $$H^*(\mc{C},F)=H^*\left(\overline{\{x\}},\epsilon_{x*}M\right)=
H^*(\{x\},M)=M.$$ 
Similarly, $H^*(\mc{C}-\mc{C}',F)=M$ if $x\not\in \mc{C}'$ and zero otherwise.
Now, $H^0(\mc{C},F)=H^0(\mc{C}-\mc{C}',F)=M$ for $x\not\in \mc{C}'$ implies the statement. 
\end{proof}
\end{lemma}

\subsection{}
We consider the case $\mc{C}=[A,V]$ and $F=\epsilon_{V!}\Lambda$. 
Lemma \ref{resolutionG} and \ref{acyclic} imply 
$$
H^*_{\{A\}}([A,V],\epsilon_{V!}\Lambda) = H^*(\Gamma_{\{A\}}G(\epsilon_{V!}\Lambda,[A,V])).
$$ 
The equalities 
$$
(\epsilon_{V!}\Lambda)_x = \begin{cases} 0 & \text{if $V\not= x$} \\ 
                                                       \Lambda &\text{if $V=x$} \end{cases} \quad 
\Gamma_{\{ A \}} \epsilon_{x*} \Lambda = \begin{cases} 0 & \text{if $A\not= x$} \\ 
                                                       \Lambda &\text{if $A=x$} \end{cases}
$$
imply 
\begin{equation} \label{comparisonDGM1}
\Gamma_{\{A\}}G(\epsilon_{V!}\Lambda,[A,V]) = \Hom(S^{\{A\},\{V\}}\otimes \Lambda,\Lambda),
\end{equation}
where $S^{\{A\},\{V\}}$ is defined in \ref{SMN}. In the language of derived categories
this means 
\begin{equation} \label{quiscohomologywithsupport}
R\Gamma_{\{A\}}([A,V],\epsilon_{V!}\Lambda) \cong \Hom(S^{\{A\},\{V\}}\otimes \Lambda,\Lambda).
\end{equation}

\subsection{} Let $\alpha:[A,V]\times [B,V] \xr{} [C,V]$ be a functor such that 
$\alpha^{-1}(V)=\{(V,V)\}$ and $\alpha^{-1}(C)=\{(A,B)\}$. Moreover, we assume 
that the extension of $\alpha$:
\begin{align*}
([A,V]\times [B,V])(m) &\xr{} [C,V](m) \\
((x_0,y_0)\xr{} \dots \xr{} (x_m,y_m)) & \mapsto 
(\alpha(x_0,y_0)\xr{} \dots &\xr{} \alpha(x_m,y_m)) 
\end{align*}
maps non-degenerated chains to non-degenerated chains. If $A,B$ intersect transversally in 
$A\cap B=C$ then the intersection $\alpha(x,y)=x\cap y$ satisfies the assumptions, only
this map $\alpha$ is important for us. 

The equality 
\begin{equation}\label{supportandpushforward}
R\Gamma_{\{(A,B)\}}([A,V]\times [B,V],\bullet) = R\Gamma_{\{C\}}([C,V],\bullet) \circ R\alpha_*
\end{equation}
can be shown by considering sheaves of the form $\epsilon_{(x,y)*}M$. More 
generally, for any continuous map $\alpha:X\xr{} Y$ and a closed set $Z\subset Y$ 
we know that 
$$
R\Gamma_{\alpha^{-1}Z} = R\Gamma_{Z} \circ R\alpha_*
$$
holds (it is sufficient to consider injective sheaves which can be checked easily).

The natural map $F\xr{} R\alpha_*\alpha^*F$ yields
\begin{equation} \label{pullbackwithsupport}
R\Gamma_{\{C\}}([C,V],F) \xr{} R\Gamma_{\{(A,B)\}}([A,V]\times [B,V],\alpha^*F).
\end{equation}

Via the quasi-isomorphism from Lemma \ref{lemmaFquisG} the  map 
$F\xr{} R\alpha_*\alpha^*F$ can be rewritten as:
\begin{equation}\label{N}
N:G(F,[C,V]) \xr{} \alpha_*G(\alpha^{*}F,[A,V]\times [B,V])
\end{equation}
with 
$$
{\rm pr}_x \circ N = (\epsilon_{\alpha(x_0)*}F_{\alpha(x_m)} \xr{} \alpha_*\epsilon_{x_0*}(\alpha^*F)_{x_m} ) \circ {\rm pr}_{\alpha(x)}
$$
for all $x$. Note that $(\alpha^*F)_{x_m}=F_{\alpha(x_m)}$, and that 
$\epsilon_{\alpha(x_0)*} \xr{} \alpha_* \epsilon_{x_0*}$ is induced by 
$\epsilon_{x_0}^*\alpha^*\epsilon_{\alpha(x_0)*}=id$. In order to check that
$N$ is the correct morphism we have to show that the diagram
$$
\xymatrix
{
\alpha^*G^0(F,[C,V]) \ar[r]^-{N}&
G^0(\alpha^*F,[A,V]\times [B,V])
\\
\alpha^*F \ar[u] \ar[ur]
}
$$
commutes, which is clear.

We may compute \ref{pullbackwithsupport} by applying $\Gamma_{\{C\}}$ to \ref{N}.
For $F=\epsilon_{V!}\Lambda$ we get the dual of the morphism
\begin{equation} \label{pullbackwithsupportexplicit}
S^{\{(A,B)\},\{(V,V)\}}\otimes \Lambda \xr{} S^{\{C\},\{V\}} \otimes \Lambda
\end{equation}
induced by $\alpha$. Explicitly it maps
$$
\sum_{x} c_x x \mapsto \sum_{x} c_x \alpha(x),
$$
where $x$ runs over all chains $([A,V]\times [B,V])(m)_{{\rm nd}}$. 

\subsection{}
Again by using resolutions $G(\bullet,\bullet \bullet)$ it is easy
to see that the K\"unneth formula 
\begin{multline} \label{Kuennethwithsupport}
R\Gamma_{\{A\}}([A,V],\epsilon_{V!}\Lambda) \otimes^{\mathbb{L}} R\Gamma_{\{B\}}([B,V],\epsilon_{V!}\Lambda) \xr{\cong} \\ R\Gamma_{\{(A,B)\}}([A,V]\times [B,V],\epsilon_{(V,V)!}\Lambda)
\end{multline}
is given by the map \ref{extcup}. To prove this, one has to lift \ref{extcupformula} 
to sheaves. This is straightforward and we omit the proof.   

\subsection{}
Let $A$ and $B$ intersect transversally in $C=A\cap B$, and let $\alpha:[A,V]\times [B,V] \xr{}
[C,V], (x,y) \mapsto x\cap y,$ be the intersection.
  
In \cite[Remarque~4.3]{DGM} it is conjectured that the cup product 
$$
h^*(A)\otimes h^*(B) \xr{\cup} h^*(C)
$$
is given (up to Tate twist) by the dual of  
\begin{align*}
R\Gamma_{\{C\}}([C,V],\epsilon_{V!}\Lambda) &\xr{\text{\ref{pullbackwithsupport}}} 
R\Gamma_{\{(A,B)\}}([A,V]\times [B,V],\epsilon_{(V,V)!}\Lambda)
 \\
&\xr{\text{\ref{Kuennethwithsupport}}} R\Gamma_{\{A\}}([A,V],\epsilon_{V!}\Lambda) \otimes^{\mathbb{L}} 
R\Gamma_{\{B\}}([B,V],\epsilon_{V!}\Lambda).
\end{align*}
Via the quasi-isomorphism \ref{quiscohomologywithsupport} this map equals 
\begin{align}\label{formulaofDGMexplicit}
S^{\{A\},\{V\}}\otimes S^{\{B\},\{V\}} \otimes \Lambda \xr{(\text{dual of \ref{extcup}})^{-1}}&  
S^{\{(A,B)\},\{(V,V)\}} \otimes \Lambda \\
\xr{(\text{\ref{pullbackwithsupportexplicit}})}& S^{\{C\},\{V\}} \otimes \Lambda. \nonumber
\end{align}
According to Corollary \ref{satz} we know  
$$
h^*(A)=H^{*-2\cd A}(S^{\{A\},\{V\}}\otimes \Lambda)\otimes \mu_n(-\cd A),
$$
and that \ref{formulaofDGMexplicit} gives the cup product is proved in 
Corollary \ref{Multiplikation}. 
So that the conjectured formula of \cite{DGM} holds for our decomposition.  

%%%%%<<<<<<<<<<<<<<<<< revision 2

\end{document}